\documentclass[aos]{imsart}
\RequirePackage{amsthm,amsmath,amsfonts,amssymb}
\RequirePackage[round,authoryear]{natbib}

\usepackage{xr}

\usepackage[utf8]{inputenc}
\usepackage{soul}
\usepackage{amsmath, amsthm, amssymb, graphicx, multicol, array, appendix, verbatim, multirow, comment, enumerate}
\usepackage[colorlinks=true, citecolor=blue]{hyperref}
\usepackage{cleveref}
\usepackage[table,xcdraw]{xcolor}
\usepackage{hhline}
\usepackage{graphicx}
\usepackage{placeins}
\usepackage{subfigure}
\usepackage{algpseudocode, algorithm}

\startlocaldefs

          \def\cI{{\cal  I}}

\def \RR	{\mathbb{R}}

\def \EE {\mathbb{E}}

\def \tilde{\widetilde}

\newtheorem{theorem}{Theorem}[section]
\newtheorem{corollary}{Corollary}[theorem]

\theoremstyle{remark}
\newtheorem{remark}{Remark}[section]

\endlocaldefs

\makeatletter

\newcommand*{\addFileDependency}[1]{
\typeout{(#1)}
\@addtofilelist{#1}
\IfFileExists{#1}{}{\typeout{No file #1.}}
}\makeatother

\newcommand*{\myexternaldocument}[1]{%
\externaldocument{#1}%
\addFileDependency{#1.tex}%
\addFileDependency{#1.aux}%
}

\myexternaldocument{aos-supp-R01}

\begin{document}

\begin{frontmatter}

\title{Statistical Inference on Latent Space Models for Network Data}

\begin{aug}

\author[]{\fnms{Jinming}~\snm{Li}\ead[label=e1]{lijinmin@umich.edu}}, 
\author[]{\fnms{Shihao}~\snm{Wu}\ead[label=e2]{wshihao@umich.edu}}, 
\author[]{\fnms{Chengyu}~\snm{Cui}\ead[label=e3]{chyc@umich.edu}}, 
\author[]{\fnms{Gongjun}~\snm{Xu}\ead[label=e4]{gongjun@umich.edu}}, 
and 
\author[]{\fnms{Ji}~\snm{Zhu}\ead[label=e5]{jizhu@umich.edu}}
\address[]{
Department of Statistics, University of Michigan, Ann Arbor \\ \printead[presep={\ }]{e1,e2,e3,e4,e5}}
\end{aug}



\begin{abstract}
    Latent space models are powerful statistical tools for modeling and understanding network data. While the importance of accounting for uncertainty in network analysis has been well recognized, the current literature predominantly focuses on point estimation and prediction, leaving the statistical inference of latent space models an open question. This work aims to fill this gap by providing a general framework to analyze the theoretical properties of the maximum likelihood estimators. In particular, we establish the uniform consistency and asymptotic distribution results for the latent space models under different edge types and link functions. Furthermore, the proposed framework enables us to generalize our results to the dependent-edge and sparse scenarios. Our theories are supported by simulation studies and further illustrated by an application involving a statistician coauthorship network. 
\end{abstract}


\begin{keyword}
\kwd{latent space models}
\kwd{maximum likelihood inference}
\kwd{network analysis}
\end{keyword}

\end{frontmatter}


\section{Introduction}

Network data have become increasingly prevalent in various fields of study and applications, such as social networks \citep{traud2012social}, biological networks \citep{bullmore2009complex}, trading networks \citep{chaney2014network}, and collaboration networks \citep{ji2016coauthorship}, among others. Such data consist of a collection of nodes, representing entities, and edges, representing the relationships or interactions between these entities. A network can be denoted by an $n \times n$ adjacency matrix $A$, where $n$ is the number of nodes and $A_{ij}$ represents the information of the edge between nodes $i$ and $j$, which can take various values depending on the specific application. For example, in a friendship network on social media, $A_{ij}$ may be binary, where $A_{ij} = 1$ if and only if two nodes are mutual friends; in neuroimaging networks, $A_{ij}$ may represent transformed brain connectivity measures
and take continuous values. In this paper, we consider the general setting where networks are undirected and have no self-loops, i.e., $A_{ij} = A_{ji}$ for all $1 \leq i < j \leq n$ and $A_{ii} = 0$ for all $1 \leq i \leq n$.

Latent space models are powerful tools to perform statistical modeling and inference on network data \citep{hoff2002latent, athreya2017statistical, smith2019geometry}. By embedding nodes in a lower-dimensional latent space, these models are able to not only capture the dependency of edges but also reveal the complex and heterogeneous network structure via the geometric relationships between the embedded positions in the latent space. One popular class of models is the inner product latent space model \citep{young2007random,zhao2012consistency,eavani2015identifying,ma2020universal,wang2023locus}, which utilizes the inner product to measure the similarity of latent positions, characterizing important network properties, such as transitivity, assortativity, and community structure prevailed in social networks \citep{smith2019geometry}, or block and crossing patterns observed in neuroimaging networks \citep{amico2018mapping}.

Many popular network latent space models fall into the following formulation. Let each node $i$ be assigned with a latent position $z_i \in \mathbb R^{r}$ in the $r$-dimensional latent space and a degree parameter $\alpha_i \in \mathbb R$ that accounts for node heterogeneity. Then for each node pair $1 \leq i < j \leq n$, we assume
\begin{equation}
    A_{ij} \sim p\left( ~ \cdot ~ | ~ \theta_{ij}\right) \mbox{ with }  \theta_{ij} = \sigma(z_i^{\top}z_j + \alpha_i + \alpha_j),
    \label{eq:general_model}
\end{equation}
where each edge $A_{ij}$ is a random variable following the distribution $p(\cdot| \theta_{ij})$ with the parameter $\theta_{ij}$ determined by the latent positions $(z_i, z_j)$ and degree parameters $(\alpha_i, \alpha_j)$ through a link function $\sigma(\cdot)$. Typically, $\sigma(\cdot)$ is a smooth and increasing function, ensuring that the model exhibits the desirable property that higher inner product similarity of latent positions and greater node heterogeneity values lead to higher expected values of $A_{ij}$.  Examples of popular models that can be formulated as Equation \eqref{eq:general_model} include:
\begin{itemize}
    \item The Random Dot Product Graph (RDPG) model \citep{young2007random, athreya2017statistical}, which considers a linear model for binary networks: $A_{ij} \overset{i.i.d.}{\sim} \text{Bernoulli}(\theta_{ij})$ with  $\sigma(x) = x$ and $\alpha = 0_{n \times 1}$. RDPG covers some popular discrete network latent space models, such as the positive definite Stochastic Block Model and its variants \citep{holland1983stochastic, zhao2012consistency}.
    \item The universal inner product 
    \citep{ma2020universal}, which incorporates the logistic link function for binary networks: $A_{ij} \overset{i.i.d.}{\sim} \text{Bernoulli}(\theta_{ij})$ with $\sigma(x) = 1 / (1 + \exp^{-x})$. This model extends the inner product model proposed in \cite{hoff2002latent} by introducing the node heterogeneity parameters ($\alpha_i, i=1,\cdots,n$).
    \item The linear low-rank decomposition type models \citep{sun2017store,wang2023locus}, which aim to model Gaussian edges such as fMRI networks: $A_{ij} \overset{i.i.d.}{\sim} \mathcal{N}(\theta_{ij}, \delta^2)$ with variance $\delta^2$, $\sigma(x) = x$ and $\alpha = 0_{n \times 1}$.
\end{itemize} 

While latent space models have been widely used for modeling network data, the important yet challenging statistical inference question remains unaddressed.  This is not only crucial for quantifying the estimation uncertainty of latent positions, but can also facilitate downstream inference tasks, such as link prediction and network testing.
In the existing literature, the asymptotic distribution is only studied for the spectral embedding estimators under the RDPG model \citep{athreya2017statistical}.
However, RDPG, which models binary edges through a linear model, may not be the most suitable choice.
\cite{ma2020universal} models binary edges with a logistic link function and establishes {\it average} consistency result of latent positions. However, deriving asymptotic distribution results for latent positions under the logistic link setting, or other link functions and edge types in general, remains a challenging problem due to the non-linearity structures. The analysis techniques developed in the RDPG literature largely rely on the linear model structure and the corresponding spectral decomposition, which are not applicable to general settings with different edge types and link functions. Moreover, the properties of latent position estimators, under more realistic settings such as edge dependency and sparsity, have yet to be investigated. 

{
Conducting inference on latent space models also poses unique technical challenges and differs fundamentally from the existing literature for parametric models with a diverging number of parameters. In particular, the identifiability issue in latent space models requires additional consideration, unlike statistical inference for generalized linear models with a diverging number of parameters \citep{fahrmeir1985consistency}. We address this by considering a constrained MLE under the identifiability conditions and introducing a novel Lagrangian-penalty strategy designed for latent space models.
Moreover, in the literature of classical exponential families \citep{portnoy1988asymptotic} and more general $M$-estimation \citep{he2000parameters}, valid inference typically requires that the ratio between the square of the number of parameters and the sample size goes to zero. Our setup, however, violates this requirement because the ratio converges to a constant, rendering the existing techniques for parametric models with a diverging number of parameters inapplicable.

}

In this work, we establish statistical inference results for the maximum likelihood estimators of a general class of network latent space models. Our main contributions are summarized as follows.

\begin{enumerate}
    \item We propose a unified and flexible framework for analyzing the theoretical properties of the maximum likelihood estimators of network latent space models under a general model setup. Our approach is based on the examination of the high-dimensional structures of the score vector and the Hessian matrix of the log-likelihood function, which is fundamentally different from the techniques used in the existing literature \citep[e.g.][]{athreya2017statistical,ma2020universal}. To overcome the theoretical challenge arising from the singularity of the Hessian matrix due to the identifiability issue, we introduce a Lagrange-type adjustment to the likelihood function to accommodate the identifiability constraints, which allows for the development of a comprehensive statistical inference framework for various network data settings. 
    
    \item  To the best of our knowledge, this work is the first to establish the uniform consistency as well as asymptotic distribution results for the maximum likelihood estimators
    under general edge types and link functions. Notably, we show that the asymptotic variance is optimal in the sense of achieving the Cramer-Rao information lower bound for individual node estimation. {Specifically, under independent-edge settings, the asymptotic covariance matrix for the estimators of node degree $\alpha_i$ and node latent position $z_i$ matches the inverse of the oracle Fisher information matrix of these parameters given that all node degree and latent position parameters associated with the other $n-1$ nodes are \emph{known}.} The established asymptotic distribution results would also pave the way for downstream inference tasks in network data applications, such as constructing confidence intervals for link predictions.

    \item Adopting the same analysis framework, we further extend our results towards more practical settings, where the network is allowed to be sparse and edge variables are allowed to have additional dependencies beyond those explained by the latent positions. {We show that our results adapt to the optimal sparse regime in network analysis up to an arbitrarily small polynomial order, i.e., $\Theta(n^{-1+\epsilon})$, under independent-edge settings. We also examine the tolerance of our estimator against edge dependency and the related asymptotic properties, where the estimation method can be viewed as a pseudo-likelihood approach.}
\end{enumerate}

The rest of this paper is organized as follows. We start by discussing the network model setting with independent edges in Section \ref{sec:ind_case}, where we present the problem setup, the proposed theoretical framework, and the main theoretical results. In Section \ref{sec:generalized_settings}, we extend the analysis towards the dependent-edge and sparse-edge settings. Comprehensive simulation results are presented in Section \ref{sec:simulation}, and a data application involving a statistician coauthorship network is presented in Section \ref{sec:real_data}. Additional numerical results and proofs are presented in the Supplementary Material.

\section{Maximum Likelihood Inference via the Lagrange-Adjusted Hessian}
\label{sec:ind_case}

\subsection{Problem Setup}
\label{sec:problem_setup_notations}

Our goal is to analyze the properties of the maximum likelihood estimators of the general latent space model introduced in Equation \eqref{eq:general_model}, where the parameters of interests are $z_i\in \mathbb R^{r}$ and $\alpha_i\in \mathbb R$ for $1 \leq i \leq n$. Following the latent space model literature \citep[e.g.][]{ma2020universal},
we treat $z_i$ and $\alpha_i$ as fixed parameters and do not pose specific distributional assumptions on how they are generated. For situations where these parameters are treated as latent random variables, our analysis could be regarded as making statistical inferences on the realized values of those latent variables.

We use the following notations throughout the analysis. Let $Z_{n \times r}$ represent the latent position matrix, with the $i$-th row $z_i^{\top}$ being the latent position vector for node $i$. Similarly, let $H_{n \times (r + 1)}$ be a matrix with the $i$-th row   defined as $h_i^{\top} = (z_i^{\top}, 1)$. Let $\alpha_{n \times 1}$ be the vector of node heterogeneity parameters. Denote $\phi$ as an $n(r + 1)$-dimensional vector, whose $i$-th block $\phi_i = (z_i^{\top}, \alpha_i)^{\top}$ is the $(r+1)$-dimensional vector that includes all latent parameters associated with node $i$. Denote $\pi_{ij} = z_i^{\top}z_j + \alpha_i + \alpha_j$ such that we can write $\theta_{ij} = \sigma(\pi_{ij})$, and in matrix form, denote $\Pi_{n \times n} = ZZ^{\top} + \alpha 1_n^{\top} + 1_n \alpha^{\top}$, where $1_n$ is the $n$-dimensional vector with all entries being $1$. Similarly, we use $0_n$ to denote the $n$-dimensional zero vector.

For each edge $A_{ij}$, $1 \leq i < j \leq n$, denote the log-likelihood function with respect to $\pi_{ij}$ as 
$$l(\pi_{ij}; A_{ij}) = \log p(A_{ij} | \theta_{ij}) = \log p(A_{ij} | \sigma(\pi_{ij})).$$
{We also define $A_{ij} = A_{ji}$, $\theta_{ij} = \theta_{ji}$, and $\pi_{ij} = \pi_{ji}$ for $j<i$, which will be used in later analysis.}
For simplicity, we write $l(\pi_{ij}; A_{ij})$  as $l(\pi_{ij})$ and further introduce $l'(\pi_{ij}), l''(\pi_{ij}), l'''(\pi_{ij})$ to denote the $1^{\mathrm{th}}, 2^{\mathrm{nd}}, 3^{\mathrm{rd}}$-order derivatives of $l(\pi_{ij}; A_{ij})$ with regard to $\pi_{ij}$, respectively. Given $\Pi_{n \times n}$, if all the edges are independent, the log-likelihood function of the network data can be written as
\begin{equation*}
    L(\phi; A) = \sum_{i  = 1}^n \sum_{j = i + 1}^n l(\pi_{ij}).
\end{equation*}  

We use $\phi^*$, $\phi_i^*$, $z_i^*$, $\alpha_i^*$,  $\pi_{ij}^*$, and $ \theta^*_{ij}$  to denote the true parameters. We use $\|\cdot\|$, $\|\cdot\|_1$, $\|\cdot\|_F$, $\|\cdot\|_{\max}$ to denote matrix 2-norm, 1-norm, Frobenius norm, and max norm, correspondingly.
Following the literature \citep[e.g.][]{ma2020universal}, we consider the maximum likelihood estimators under boundedness constraints. Specifically, for some large constant $M$,   we seek to maximize $  L(\phi; A)$ subject to 
$\sup_{1 \leq i \leq n} \|\phi_i\| \leq M$.
{To address the identifiability issue of the model parameters, we impose additional constraints on the estimators such that $Z$ is centered with $Z^{\top}1_{n}=0_r$ and $Z^{\top}Z/n$ is diagonal \citep{zhang2022joint}. The rationale for these constraints is explained in the next section. In the following, we use $\hat\phi$, $\hat\phi_i$, $\hat z_i$, $\hat \alpha_i$,  $\hat \pi_{ij}$, and  $\hat \theta_{ij}$  to denote the maximum likelihood estimators that maximize $L(\phi; A)$ under the boundedness and identifiability constraints. 
} 


\subsection{Assumptions}\label{sec:assum}
For the purpose of theoretical analysis, we make the following assumptions. 

\begin{enumerate}[I.]
    \item All edge variables $A_{ij}$, $1 \leq i < j \leq n$, are independent conditioning on $\phi^*$.
    \smallskip
    \item There exists a large  $M > 0$ such that the parameter set $\Xi = \left\{\phi \mid \sup_{1 \leq i \leq n} \|\phi_i\| \leq M \right\}$ contains the true parameters $\phi^*$. 
  \smallskip
    \item There exists a positive definite $r \times r$ matrix $\Sigma_Z$ such that ${Z^*}^{\top}Z^*/n \to \Sigma_Z$ as $n \to \infty$. {Moreover, for any $i\in[n]$, there exists a positive definite $r\times r$ matrix $\Sigma_i$ such that $ - n^{-1}\sum_{j:j\ne i} l''(\pi_{ij}^*)h_j^*h_j^{*\top} \overset{p}{\to}\Sigma_i$ as $n\to\infty$.} 
    \smallskip
    \item $Z^*$ is centered such that ${Z^*}^{\top}1_n = 0_r$, and the limiting covariance matrix $\Sigma_Z$ is diagonal with unique eigenvalues.  
    \smallskip    
    \item $l'''(\cdot)$ exists within $\Xi$. Furthermore, there exists $0< b_L < b_U$ such that $b_L \leq -l''(\cdot) \leq b_U$ and $|l'''(\cdot)| \leq b_U$ within $\Xi$.
    \smallskip
    \item There exist $t > 0$ and $s > 0$ such that for all $x \geq 0$, $P(|l'(\pi_{ij}^*)| > x ) \leq \exp(-(x/t)^s)$.

\end{enumerate}

Assumption I is a common independence assumption made by the existing literature \citep{ma2020universal,young2007random,hoff2002latent}. Assumption II is on the boundedness of the true parameters for the purpose of theoretical analysis and also commonly assumed \citep{ma2020universal}. In Section \ref{sec:generalized_settings}, we will discuss how to relax these two assumptions. Assumption III poses regular conditions on the asymptotic behavior of the {second moment} of latent positions {$\{z_i^*\}_{i\in[n]}$}. {The second part of Assumption~III is introduced to handle certain marginal cases of link functions, thereby ensuring the generality of our theoretical results. For commonly used link functions such as Gaussian, logistic, and Poisson, this part of the assumption is not needed.
}



{Assumption IV poses identifiability conditions that enable statistical inference for individual parameters in $\phi$:
(i) $Z^*$ needs to be centered, as otherwise $\alpha^*$ is not identifiable; (ii) the diagonality requirement of $\Sigma_Z$ and the uniqueness of eigenvalues avoid the orthogonal transformation issue of $Z^*$.
Specifically, 
consider any $\mu_{Z}\in\RR^r$ and orthogonal matrix $Q\in\RR^{r\times r}$. If we define 
\begin{equation} \label{eq:trans}
z_i' = Q(z_i^*-\mu_{Z}) \text{ and } \alpha_i' = \alpha_i^* + \mu_{Z}^{\top}z_i^* - \mu_Z^{\top}\mu_{Z}/2 \text{ for all } i\in[n],
\end{equation}
then $z_i'^{\top}z_j' + \alpha_i' + \alpha_j' = {z_i^*}^{\top}z_j^* + \alpha_i^* + \alpha_j^*$ for all $1\le i<j\le n$. Thus, $\alpha_i^*$'s are identifiable up to a shift $\mu_Z$ in $z_i^*$'s, and the latent positions $z_i^*$'s are identifiable up to an orthogonal transformation $Q$, because these transformations do not change the likelihood value. 
In Assumption IV, point (i) addresses the shift issue by fixing $\mu_Z$ and point (ii) fixes the rotation $Q$ via the uniqueness of eigen-decomposition. We want to emphasize that other identifiability conditions guaranteeing the entry-wise identifiability and estimability of \(Z\) and \(\alpha\) can also be adopted within our framework. For example, the diagonality assumption may be replaced by some other \(r(r-1)/2\) rotation constraints; hence, our results do not rely on the columns of \(Z^*\) being linearly independent. We use the identifiabilty conditions in Assumption IV  for their popularity in the literature \citep{ma2020universal,zhang2022joint,li2025high}, 
and for the analysis convenience. For any other identifiability
condition, our results still hold up to the corresponding transformation. Specifically, if the true \(Z^*\) and \(\alpha^*\) do not satisfy the centering and diagonality conditions in Assumption~IV, one can apply a unique transformation to \(Z^*\) and \(\alpha^*\) so that these conditions are satisfied, and all of our theoretical results then hold for the transformed true parameters.

Moreover, when the eigenvalues of \(n^{-1}Z^{*\top}{Z^*}\) are not unique (i.e., there exists an eigensubspace of dimension \(r' > 1\) that shares the same eigenvalue), additional constraints are needed to fix the rotation. Specifically, we require \(r'(r'-1)/2\) extra constraints to ensure the identifiability of the entries in the corresponding columns of \(Z^*\). 
For the main results of this paper, we focus on the identifiable cases with unique eigenvalues.
In scenarios where the eigenvalue uniqueness assumption does not hold, we show that the maximum likelihood estimator still achieves consistent eigenspace estimation, and there exists an orthogonal transformation on the true latent vectors such that our inference results remain valid. See Section~\ref{sec:bounded_case} for further details.
Finally, we would like to highlight that Assumption~IV is mainly needed for statistical inference on the individual \(Z\) and \(\alpha\) parameters due to their identification requirement. If our inference target is the \(\theta_{ij}\) parameters instead, Assumption~IV is no longer needed and all our asymptotic results on the estimator of \(\theta_{ij}\) still hold. 
}

{Assumptions V-VI concern the properties of the link function. }
Assumption V requires the link function $\sigma(\cdot)$ to be smooth and the log-likelihood function $l(\pi_{ij})$ is concave with respect to $\pi_{ij}$, which is satisfied by many widely used link functions such as {logistic}, Poisson, and linear links. Assumption VI ensures that the tail density of the edge variables decays exponentially, and similar assumptions are often made in the literature \citep{bai2013statistical}.  
The following remark provides an example under the logistic link function. 


{
\begin{remark}[Assumptions V-VI under the logistic link function.]\label{rem:log1}
    Under a latent space model with the logistic link, $A_{ij}$ is Bernoulli with mean $1/(1+\exp(-\pi_{ij}^*))$. 
    Note that $l'(\pi_{ij}^*) = A_{ij} - 1/(1+\exp(-\pi_{ij}^*))$, $l''(\pi_{ij}^*) = -\exp(\pi_{ij}^*) / (1+ \exp(\pi_{ij}^*))^2$, and $l'''(\pi_{ij}^*) = \exp(\pi_{ij}^*)(1-\exp(\pi_{ij}^*))/(1+\exp(\pi_{ij}^*))^3$.     
Then when $\Xi$ is a bounded set,  Assumption V is satisfied. 
Note that $l'(\pi_{ij})$'s are Bernoulli random variables. Assumption VI is also naturally satisfied. 
\end{remark}
}

\subsection{Theoretical Analysis via  Lagrange-Adjusted Hessian}

\label{sec:bounded_case}

In this section, we illustrate the roadmap of our analysis on the properties of the maximum likelihood estimators while leaving the detailed proofs in the Supplementary Material. Before we dive into the details, it is worthwhile to first understand the structure of the Hessian matrix of the log-likelihood function, $H(\phi) = \frac{\partial^2 L}{\partial \phi \partial \phi^{\top}}$, as it plays the central role in our analysis. Following the block structure of $\phi$,  we can present the  $n(r + 1)\times n(r + 1)$ Hessian matrix $H(\phi)$  as $n\times n$ matrix blocks with size  $(r + 1) \times (r+1)$, i.e.,  $ \left[H(\phi)\right]_{ij}  = \frac{\partial^2 L}{\partial \phi_i \partial \phi_j^{\top}}$ for $1\leq i, j\leq n$. 
We can further write $H(\phi)$ into two parts such that $H(\phi) = H_L(\phi) + J_L(\phi)$, where each $(r + 1) \times (r+1)$ block is:
\begin{equation} \label{eq:Hessian_H_L_J_L}
    \begin{aligned}
    \left[H_L(\phi)\right]_{ii} &= \sum_{j: j \neq i} l''(\pi_{ij}) h_j h_j^{\top}, \quad \left[J_L(\phi)\right]_{ii} = 0_{(r + 1) \times (r + 1)}, \quad 1 \leq i \leq n; \\
    \left[H_L(\phi)\right]_{ij} &= l''(\pi_{ij}) h_j h_i^{\top}, \quad \left[J_L(\phi)\right]_{ij} = l'(\pi_{ij}) \begin{pmatrix}
        I_r & 0_r \\ 0_r^{\top} & 0
    \end{pmatrix},
    \quad 1 \leq i ,j \leq n.
    \end{aligned}
\end{equation}
{For convenience, in the following, we use similar notations  $[\cdot]_i$ and  $[\cdot]_{ij}$  to denote the $i$-th block of a vector corresponding to $\phi_i$ and the $(i,j)$-th block of a matrix corresponding to ($\phi_i$, $\phi_j$).}
To establish asymptotic properties of the maximum likelihood estimators, our main idea is to first prove the first-order condition $S(\hat{\phi}) = 0$ and utilize:
\begin{equation}
    0 = S(\hat \phi) = S(\phi^*) + H(\phi^*) \times (\hat \phi - \phi^*) + R(\hat \phi, \phi^*),
\label{eq:first_order_condition}
\end{equation}
where we recall that $S(\phi) = \frac{\partial L}{\partial \phi}$ is the score vector and $H(\phi)$ is the Hessian matrix of the log-likelihood function, and $R(\hat \phi, \phi^*)$ denotes the residual term that depends on both $\hat \phi$ and $\phi^*$. Ideally, we would want the following {three} conjectures to hold. {Conjecture (1):} the residual term would be uniformly of smaller order and thus be negligible; {Conjecture (2):} the Hessian matrix would be negative definite after proper scaling; {and Conjecture (3):  $[H^{-1}(\phi^*)S(\phi^*)]_i\approx [\{H_L(\phi^*)\}_{ii}]^{-1}[S(\phi^*)]_i$.}
If these conjectures were true, we could show that for each $1 \leq i \leq n$, 
\begin{align*}
    \sqrt{n} \big(\hat \phi_i - \phi_i^*\big) &= -\sqrt{n} \left[H^{-1}(\phi^*)S(\phi^*)\right]_i + o_p(1) \\ &
  = \left[-n^{-1} \{H_L(\phi^*)\}_{ii}\right]^{-1}[n^{-1/2}S(\phi^*)]_i + o_p(1),
\end{align*}
where 
the last formula would further allow us to establish the individual asymptotic distribution using the central limit theorem. 

However, there are major theoretical challenges in proving the above conjectures. For Conjecture (1), establishing asymptotic distribution for every $1 \leq i \leq n$ requires a uniform upper bound of $[R(\hat \phi, \phi^*)]_i$, which means that we need to go beyond the average consistency of $\hat \phi$ and prove uniform consistency. We overcome this challenge by considering the integral form mean value theorem for vector-valued functions, 
which can be seen as a lower-order expansion of the first-order condition in Equation \eqref{eq:first_order_condition}:
\begin{equation}
    0 = S(\hat \phi) = S(\phi^*) + \tilde H \times (\hat \phi - \phi^*),
    \label{eq:first_order_condition_lower_order}
\end{equation}
where $\tilde H = \int_0^1 H(\phi^* + t(\hat \phi - \phi^*))dt$. Thus, if we were able to show that the Hessian matrix is negative definite after proper rescaling, we could characterize the uniform norm of $\hat \phi - \phi^*$. Hence, the second challenge that needs to be solved is to verify Conjecture (2) within a small neighborhood of $\phi^*$. 

Unfortunately, it turns out that the negative definite statement is not true. As shown in the proofs in the Supplementary Material, the leading term of the Hessian matrix, $H_L(\phi^*)$, has exactly $r(r + 1) / 2$ eigenvalues equal to $0$, which coincides with the number of constraints required for the model identifiability. 
To address similar challenges caused by singular Hessian matrices, 
the idea of introducing auxiliary Lagrange multiplier terms has been adopted in the literature of maximum likelihood inference. For example, \cite{aitchison1958maximum} and \cite{silvey1959lagrangian}  utilized the idea to study the constrained maximum likelihood inference problems with low-dimensional models. \cite{el1994wald} further explored the technique of using Lagrangian multiplier terms to study the likelihood functions with singular information matrices, while still limiting the discussion in the classic low-dimensional settings. Recently, \cite{wang2022maximum} adopted a similar idea to analyze the maximum likelihood estimators of the high-dimensional factor analysis models. We want to emphasize that compared with the existing works, analyzing network latent space models has several unique challenges that are not trivial to overcome. For instance, the latent positions of network latent space models jointly appear as $ZZ^{\top}$, which causes significantly more challenges in studying the structure of the Hessian matrix, especially the characterization of its matrix 1-norm. Moreover, due to the distinctive properties of network data, our analysis involves models with node heterogeneity parameters and considers the sparse network setting (see details in Section \ref{sec:generalized_settings}), which introduces additional challenges in the theoretical analysis that the existing techniques cannot be applied to directly.

Motivated by the existing literature, for our theoretical analysis, we consider a modified loss function that combines the log-likelihood function with additional Lagrange-type penalty terms  corresponding to the identifiability constraints of the network model. Formally, we consider the loss function below:
\begin{equation}
    Q(\phi) = L(\phi) + P(\phi) := \sum_{i = 1}^n \sum_{j = i + 1}^n l(\pi_{ij}) - \frac{c}{2}n^2 \left\|\text{upper\_diag}(Z^{\top}Z/n)\right\|_F^2 - \frac{c}{2}n^2\Big\|\frac{Z^{\top}1_n}{n}\Big\|_F^2,
    \label{eq:loss_function_bounded}
\end{equation}
where $L(\phi) = \sum_{i = 1}^n \sum_{j = i + 1}^n l(\pi_{ij}) $ is the log-likelihood function, $P(\phi)$, which is defined as the last two terms in Equation \eqref{eq:loss_function_bounded}, poses penalties on the off-diagonal entries of $Z^{\top}Z/n$ as well as the mean vector of $Z$, and $c$ is a  multiplier satisfying $0 < c < b_L$, with $b_L$ specified in Assumption V. The two penalty terms in $P(\phi)$ correspond to the rotation and centering identifiability conditions listed in Assumption IV. In the case where different rotation identifiability conditions are used, our analysis can be adapted accordingly by changing the first penalty term in $P(\phi)$. 

We would like to highlight that the estimators of the model parameters, denoted by $\hat\phi$, obtained by maximizing $Q(\phi)$ are the corresponding maximum likelihood estimators introduced in Section \ref{sec:problem_setup_notations} that satisfy the desired identifiability constraints.
In particular, the task of maximizing  Equation \eqref{eq:loss_function_bounded} is equivalent to the following two-step procedure:
  first obtaining a set of estimators of the model parameters from maximizing $L(\phi)$ only, and then transforming the obtained set of estimators (particularly $\hat Z$) such that they satisfy the considered identifiability conditions.
   This is because, due to the identifiability issue,  within the equivalent class of all sets of ``nonidentifiable'' estimators, they all yield the same value of $L(\phi)$, while the identifiability constraints further refine them to the set of estimators corresponding to directly optimizing Equation \eqref{eq:loss_function_bounded}. 
Therefore, the choice of the positive constant $c$ does not affect the estimation, and practically, when calculating $\hat\phi$, we use the equivalent two-step procedure, which does not rely on $c$. 

For the same reason, we further emphasize that the additional $P(\phi)$ in Equation \eqref{eq:loss_function_bounded} does not work the same way as the penalty terms in the literature of shrinkage regression, and similarly the parameter $c$ does not function as a tuning parameter but rather a theoretical tool for us to study the properties of the maximum likelihood estimators. In particular, in the theoretical analysis, it can be proved that the Hessian matrix of $P(\phi)$ can be written as $r(r+ 1)/2$ linearly independent vectors that form a basis of the null space of $H_L(\phi)$. Thus, 
the Hessian matrix of $Q(\phi)$ is negative definite with the information provided by $P(\phi)$ and all the analysis aforementioned can be applied to the maximizers of $Q(\phi)$, as they also satisfy the first-order condition (Equation \eqref{eq:first_order_condition}).

{Finally, proving Conjecture~(3) presents unique challenges in latent space network models. To establish the approximation in Conjecture~(3), we require a careful characterization of $H^{-1}(\phi^*)$, which is the inverse of a high-dimensional matrix whose dimension grows with $n$. A particular challenge is that the diagonal blocks $[H_L(\phi^*)]_{ii}$ and the off-diagonal blocks \( [H_L(\phi^*)]_{ij} \) contribute terms of the same order of magnitude to the Hessian matrix, yet cannot be easily separated in the inverse. This challenge is unique from those in high-dimensional latent factor models \citep{bai2003inferential,wang2022maximum} with a four-block-structured Hessian, rendering the Schur-complement-type techniques used therein inapplicable. To address this, we introduce novel matrix decomposition techniques specifically for network latent space models that disentangle the diagonal and off-diagonal blocks in their respective contributions to the expansion.
Furthermore, the order of $J_L(\phi^*)$ plays a critical role in the analysis, and its asymptotic behavior is particularly complex, especially in the dependent and sparse setting discussed in Section~\ref{sec:generalized_settings}. To address these challenges, we carefully characterize moment-based quantities of $J_L(\phi^*)$ to account for the underlying dependencies, and leverage tools from sparse matrix analysis to derive theoretical insights into how sparsity impacts statistical inference under general setups.
We refer interested readers to Lemmas~B.11, B.12, and B.13 in the Supplementary Material.
}


With the proposed analysis techniques,
we can establish the following theoretical results.

{
\begin{theorem}[Uniform and Average Consistency]
\label{thm:uni_consistency_bounded}  
Under Assumptions I-VI, for any $\epsilon > 0$, we have
$$\|\hat \phi - \phi^*\|_{\max} = O_p(n^{-1/2 + \epsilon})  ~~\text{and}~~n^{-1}\|\hat \phi - \phi^*\|^2 = O_p(n^{-1}) .$$
\end{theorem}
}

\begin{remark}\label{remark:unif}
    The uniform consistency rate established in Theorem \ref{thm:uni_consistency_bounded} is nearly optimal, approaching the best rate of individual consistency $O_p(n^{-1/2})$.
    This result has not been established in the existing literature on network latent space models except under the RDPG model \citep{athreya2017statistical}. However, the techniques used to analyze RDPG rely on the linear model structure, whereas our analytical framework is applicable to the general class of latent space models. As for average consistency, the rate established in Theorem \ref{thm:uni_consistency_bounded} is optimal and aligns with the rate obtained in \cite{ma2020universal}, which focuses on binary networks with logistic link function and Bernoulli edges. However, our proof is different from \cite{ma2020universal} and utilizes convexity arguments that not only suit broader settings but also lead to new results on the uniform consistency as well as asymptotic normality property as shown in the following theorem. 
\end{remark}

{
\begin{remark}\label{remark:spectra}
    Theorem \ref{thm:uni_consistency_bounded} concerns the uniform consistency and average consistency in the estimation of $\phi^*$. As discussed before, when the eigenvalue uniqueness assumption in Assumption IV is not satisfied, i.e., when there exists an eigen subspace of $n^{-1}Z^*{Z^*}^{\top}$ with dimension greater than $1$ that shares the same eigenvalue, the entries in $\phi^*$ the correspond to those columns in $Z^*$ are not identifiable. In such cases, we show that $\{\hat{z}_i\}_{i\in[n]}$ achieves consistent eigenspace estimation. Specifically, let $\hat{U},U^*\in\RR^{n\times r}$ be the left-singular vectors of $\hat{Z}$ and $Z^*$, respectively, arranged with corresponding singular values in a non-decreasing order. Then there exists an orthogonal matrix $O\in\RR^{r\times r}$ such that $\|\hat{U} - U^*O\|_F = O_p(n^{-{1}/{2}}\log^{{1}/{2}}n)$. Moreover,  let $\cI^{\mathrm{sub}}=\{t,t+1,\ldots,s-1,s\}$ be any interval subset lying in $[r]$ for some $1\le t\le s \le r$. Let $\{\lambda_i\}_{i\in[r]}$ be the singular values of $Z^*$ in a non-decreasing order. Define $\lambda_{0} = +\infty$ and $\lambda_{r+1} = -\infty$. As long as $\lambda_{t-1}$ is distinct from $\lambda_{t}$ and $\lambda_{s}$ is distinct from $\lambda_{s+1}$, we have $\|\hat{U}_{\cI^{\mathrm{sub}}}-{U}_{\cI^{\mathrm{sub}}}^*O'\|_F = O_p (n^{-{1}/{2}} \log^{{1}/{2}}n )$ for some orthogonal matrix $O'\in\RR^{|\cI^{\mathrm{sub}}|\times |\cI^{\mathrm{sub}}|}$. These results follow from a Kahan-Davis Theorem argument \citep{yu2015useful} and our intermediate finding that $\|n^{-1}\hat{Z}\hat{Z}^{\top} - n^{-1}Z^*{Z^*}^{\top}\|_F = O_p(n^{-{1}/{2}} \log^{{1}/{2}}n )$. Please see more details in the proof at the end of Section B.4 in the Supplementary Material.
\end{remark}
}

Next, we develop theories to enhance statistical inference involving either individual or multiple nodes, a topic scarcely explored in the existing literature. Given $m$ node indices $\mathcal{I} = (i_1, i_2, ..., i_m)$, denote the following consistent estimators of {information matrices}: 
\begin{equation}
    \label{eq:variance_estimators}
        \hat\Sigma_{\mathcal{I}} = \text{diag}(\hat \Sigma_{i_1}, \hat\Sigma_{i_2},...,\hat\Sigma_{i_m}) ~ \text{with} ~
        \hat \Sigma_{i} = - \frac{1}{n}\sum_{j: j \neq i} l''(\hat \pi_{ij}) \hat h_j \hat h_j^{\top}, ~ i \in \mathcal{I}; 
\end{equation}
\begin{theorem}[Joint Asymptotic Distribution]
\label{thm:joint_distribution_bounded}
Given $m$ node indices $\mathcal{I} = (i_1, i_2, ..., i_m)$, denote $\phi_{\mathcal{I}} = (\phi_{i_1}^{\top}, \phi_{i_2}^{\top},...,\phi_{i_m}^{\top})^{\top}$. Under Assumptions I-VI,  we have
$$\sqrt n \Big[\widehat{\text{Var}(\hat \phi_{\mathcal{I}})}\Big]^{-1/2} (\hat \phi_{\mathcal{I}} - \phi_{\mathcal{I}}^*) \overset{d}{\to} \mathcal{N} (0, I_{m(r+1)}),$$
where {$\widehat{\text{Var}(\hat \phi_\mathcal{I})} = \hat\Sigma_\mathcal{I}^{-1}$ and $\hat \Sigma_\mathcal{I}$ is defined in Equation \eqref{eq:variance_estimators}.}
\end{theorem}

\begin{remark}
    {
As discussed in Section \ref{sec:assum}, the identifiability of $\phi^*$ is crucial in Theorem \ref{thm:joint_distribution_bounded} because our inference targets individual entries in the latent space. The identifiability of \(\phi^*\) is based on the identifiability conditions in Assumption~IV. When the centering and diagonality constraints are not initially met, a unique transformation can be applied to the true parameters in \(\phi^*\) to satisfy these conditions and meanwhile to ensure identifiability. Notably, similar inference results can be obtained under other identifiability conditions that guarantee the identifiability and estimability of parameters, using a similar Lagrangian-adjusted Hessian approach employed in this paper.
When the eigenvalues of \(\Sigma_Z\) are not unique, our inference results remain valid up to an unknown orthogonal transformation of the latent space. Specifically, there exists an orthogonal matrix \(O \in \mathbb{R}^{r\times r}\) such that the corresponding entries of \(Z^*\) in \(\phi^*\) are replaced by \(Z^*O\) in Theorem \ref{thm:joint_distribution_bounded}.
}
\end{remark}

In the special case of $m=1$, we have the following individual asymptotic distribution result.



\begin{corollary}[Individual Asymptotic Distribution]
\label{thm:ind_distribution_bounded}
Under Assumptions I-VI,  for each $1 \leq i \leq n$,  
$$\sqrt{n} \Big[\widehat{\text{Var}(\hat \phi_i)}\Big]^{-1/2} (\hat \phi_i - \phi_i^*) \overset{d}{\to} \mathcal{N} (0, I_{r+1}),$$
where {$\widehat{\text{Var}(\hat \phi_i)} = \hat\Sigma_i^{-1}$, with 
 $\hat \Sigma_i$  defined in Equation \eqref{eq:variance_estimators}. }
\end{corollary}

{
   The asymptotic distributions of maximum likelihood estimators derived in Corollary \ref{thm:ind_distribution_bounded} are optimal in the sense of achieving the Cramer-Rao information lower bound. To illustrate this, under a latent space model with the logistic link, and assuming all the embeddings for nodes in $[n] \setminus {i}$ are known, this asymptotic variance (more specifically, its version without plug-in estimators) matches the inverse of the Fisher information matrix for running logistic regression with $\{(z_j, A_{ij})\}_{j \in [n] \setminus {i}}$ pairs. This is also true for other commonly used link functions. Moreover, under the RDPG model, our asymptotic variance coincides with the result from the Adjacency Spectral Embedding estimator \citep{athreya2017statistical}. Nevertheless, as discussed in Remark \ref{remark:unif}, our analysis adopts a different approach, and our result holds for a more general class of network latent space models. Furthermore, the asymptotic variance can be estimated from samples, allowing for straightforward inference on $\phi_i$, such as constructing confidence regions.
}

The joint distribution result in Theorem \ref{thm:joint_distribution_bounded} is useful for various downstream inference problems, such as link prediction inference, network testing problems, and network-assisted regression. For instance, for link prediction inference with Bernoulli edge variables and logistic link function, the following result can be naturally established by utilizing Theorem \ref{thm:joint_distribution_bounded} for the special case of two nodes and Slutsky's theorem. Specifically, given $\mathcal{I} = \{i, j\}$, $1 \leq i < j \leq n$, under Assumptions I-VI, if $A_{ij}$ are Bernoulli random variables and $\sigma(x) = 1 / (1 + e^{-x})$, we have
\begin{equation}
    \label{cor:link_prediction_bounded}
    \sqrt{n} \Big[\widehat{\text{Var}(\hat \theta_{ij})}\Big]^{-1/2} (\hat \theta_{ij} - \theta_{ij}^*) \overset{d}{\to} \mathcal N \left(0, 1\right),
\end{equation}
where in this case $\theta_{ij}$ is the linkage probability between node $i$ and node $j$, 
and 
$$\widehat{\text{Var}(\hat \theta_{ij})} = \left[\hat\theta_{ij}(1-\hat\theta_{ij})\right]^2\begin{pmatrix} \hat h_j \\ \hat h_i  \end{pmatrix}^{\top} \begin{pmatrix}\hat\Sigma_i^{-1}  & 0 \\0 & \hat\Sigma^{-1}_j\end{pmatrix} \hat \Omega_{\mathcal{I}} \begin{pmatrix}\hat\Sigma_i^{-1}  & 0 \\0 & \hat\Sigma^{-1}_j\end{pmatrix}\begin{pmatrix} \hat h_j \\ \hat h_i \end{pmatrix}.$$
{
Note that the identifiability issue related to Assumption~IV  does not arise as a concern for the \(\theta_{ij}\) parameters. Rather, the identifiability conditions in the latent space provide analytical tools that facilitate inference on the latent vectors and associated parameters, thereby enabling inference on \(\theta_{ij}\). Consequently, the identifiability conditions in Assumption~IV serve only as a technical ``working'' condition and can be removed if a researcher is interested solely in \(\theta_{ij}^*\), rather than the latent space parameters.
}

{
\begin{remark}
There is a long line of work on conducting inference for maximum likelihood estimators (or more generally, \(M\)-estimators) with a diverging number of parameters. We want to emphasize some fundamental differences between our work and the literature. 
(1) First, parameters in latent space models are not directly identifiable nor estimable. One cannot simply analyze the maximizer of the log-likelihood function; in fact, it is not even unique. To address this, we impose identifiability conditions to resolve the identifiability issue and introduce the Lagrangian penalty term to manage the singularity in the log-likelihood Hessian. This challenge, along with our subsequent construction and analysis, sets our approach apart from inference on MLEs in setups like generalized linear models \citep{fahrmeir1985consistency}, even though interactions between latent vectors appear linear.
(2) Second, the growth rate of the number of parameters in latent space models lies outside the scope of traditional parametric inference for maximum likelihood estimators with a diverging number of parameters. Let \(n^*\) and \(p^*\) denote the sample size and the number of parameters in the MLE procedure. \cite{portnoy1988asymptotic} studied the asymptotic behavior of likelihood methods for exponential families when the number of parameters tends to infinity, showing normal approximation results under
$
{{p^*}^2}/{n^*} \to 0.
$
They further demonstrated by example that \({{p^*}^2}/{n^*} \to 0\) is necessary for normal approximation. Similar rate requirements appear in \cite{he2000parameters} for \(M\)-estimators in general parametric models. 
In our setup, we have \(n^* = n(n-1)/2\) and \(p^* = n(r+1)\), yielding
$
{{p^*}^2}/{n^*} \;\to\; 2(r+1)^2 \;>\; 0,
$
which does not converge to \(0\). This distinguishes our setup from statistical inference on maximum likelihood estimators in the classic parametric regime with a diverging number of parameters and highlights its non-triviality.  (3) Beyond the difficulties of analyzing the MLE in latent space models, our work contributes to fill the gap in the literature regarding inference for latent space models of network data with general edge types. Furthermore, we extend our framework to sparse and dependent-edge settings in the next section, thereby broadening the scope of statistical inference on latent space models for network data.
\end{remark}
}



\section{Generalization towards Dependent-Edge and Sparse Settings}
\label{sec:generalized_settings}

In this section, we further extend our analysis towards the dependent-edge and sparse network settings, which represent more realistic scenarios but also pose greater theoretical challenges. To the best of our knowledge, we are the first to establish the asymptotic properties of the maximum likelihood estimators under these settings.

\subsection{Results for Dependent-Edge Networks}\label{sec:dept_settings_theory}

In the existing literature on network latent space models, a common assumption being made is that the latent positions are sufficient to account for the dependency structure among edges. Thus, edge variables $A_{ij}$ are usually assumed to be conditionally independent as in Assumption I. However, this assumption in practice does not always hold for reasons such as mis-specifying the dimension of latent space. We therefore consider the setting where edge variables $A_{ij}$ are allowed to have additional weak dependencies after conditioning on $\phi$. Under these settings, the network latent space model can be seen as a proxy working model and the estimators are the maximizers of the pseudo-likelihood function. With the same analysis techniques on the Lagrange-adjusted Hessian, we demonstrate that theoretical results can still be established. 

For the theoretical analysis, we relax Assumption I and substitute it with the following assumptions related to the boundedness of a few moment quantities and the convergence of score subvectors. These assumptions restrict the edge dependencies from being overly strong, which could otherwise influence the validity of inference.

\begin{enumerate}[I.]
    \setcounter{enumi}{6}
    
    \item There exists $M > 0$ large enough such that
    \begin{enumerate}[i.]
    \smallskip
        \item $n^{-2}\sum_{i=1}^n \sum_{j=1}^n \big\{\frac{1}{\sqrt{n}}\sum_{m:m\ne i,j}\mathbb{E}(l'(\pi_{mi}^*) l'(\pi_{mj}^*))\big\}^2 \leq M$.
        \smallskip
        \item $n^{-2}\sum_{i=1}^n\sum_{j=1}^n\mathbb{E}\big[\frac{1}{\sqrt n}\sum_{m:m\ne i,j}\big\{l'(\pi_{mi}^*)l'(\pi_{mj}^*) -\mathbb{E}(l'(\pi_{mi}^*)l'(\pi_{mj}^*))\big\}\big]^2 \leq M$.
      
    \end{enumerate}
    \smallskip
    \item There exists $M > 0$ large enough and $\zeta > 2$ such that
    
    \smallskip
    $\mathbb{E}({n}^{-1}\sum_{i=1}^n\|\frac{1}{\sqrt{n}}\sum_{j: j \neq i} l'(\pi_{ij}^*)h^*_j\|^{\zeta}) \leq M$. 
    \smallskip
    \item There exists $M > 0$ large enough such that for all $1 \leq s \leq n$ and $1\le q\le r$:
    \begin{enumerate}[i.]
        \smallskip
        \item $\EE\|{n}^{-1}\sum_{j:j\ne s}\sum_{m: m \neq j} (n^{-1}\sum_{l:l\ne j}l''(\pi_{jl}^*)h^*_l{h^*_l}^{\top})^{-1} l'(\pi_{jm}^*)h^*_m{h^*_j}^{\top} l''(\pi_{js}^*)\|^2 \leq M$,
        \smallskip
        \item $\EE\|{n}^{-1}\sum_{j: j \neq s}\sum_{m: m \neq j} (n^{-1}\sum_{l:l\ne j} l''(\pi_{jl}^*)h^*_lh^{*\top}_l)^{-1} l'(\pi_{sj}^*) l'(\pi_{jm}^*)h^*_m\|^2 \leq M$, 
        \smallskip
        {
        \item $\EE\big\|n^{-1}\sum_{j=1}^n\sum_{m: m \neq j}h_{jq}^*\big(n^{-1}\sum_{l: l\neq j} l''(\pi_{jl}^*)h_l^*{h_l}^{*\top}\big)^{-1}l'(\pi_{jm}^*)h_m^*\big\|\le M$.
        }
    \end{enumerate}
     \smallskip
    {\item  For an integer $m$ and any $m$ node indices $\mathcal{I} = (i_1, i_2, ..., i_m)$, there exists a non-random sequence $\Omega_{\mathcal{I}}:=\Omega_{n,\cI}$ such that $n^{-1/2}\Omega_{\mathcal{I}}^{-1/2}\big([S(\phi^*)]_{i_1}, [S(\phi^*)]_{i_2}, ..., [S(\phi^*)]_{i_m} \big) \overset{d}{\to} \mathcal{N}(0, \\I_{m(r+1)})$ as $n\to\infty$, where  $S(\phi^*) = \frac{\partial L}{\partial \phi}|_{\phi=\phi^*}$ is the score vector evaluated at the true parameters $\phi^*$, and $[S(\phi^*)]_{i}$ denotes the subvector of  $S(\phi^*)$   corresponding to  all latent parameters associated with node $i$, i.e., $\phi^*_i$.}
\end{enumerate}

Assumptions VII-IX aim to characterize the strength of dependency among edge variables. Specifically, Assumptions VII and VIII restrict the dependent behavior of the score vector and Hessian matrix and are useful to establish uniform consistency results, whereas Assumption IX limits the order of residual term, $R(\hat \phi, \phi^*)$, in Equation \eqref{eq:first_order_condition} and hence leads to the asymptotic distribution results under the dependent-edge setting. {Assumption X assumes the asymptotic normality of a finite collection of subvectors of $S(\phi^*)$. It resembles common low-dimensional score convergence assumptions in the factor analysis literature \citep{bai2003inferential,wang2022maximum}, which serve to characterize score convergence with additional dependence without specifying the exact distribution form of that dependence.}

{
Under the independent-edge setup in Section \ref{sec:bounded_case}, Assumption X holds by standard central limit theorems.
  Assumptions~VII--IX also hold for any finite $\zeta>2$ because all $l'(\pi_{ij}^*)$ terms are independent and have bounded moments and exponentially decaying tails due to Assumption VI. In the dependent-edge setting, Assumption~VII concerns the scaled average second moments among the $l'(\pi_{ij}^*)$ terms and the convergence of these scaled average second moments to their population counterparts. Assumption~VIII, noting the boundedness of $h_i^*$, is analogous to Assumption~VII(i) with $\zeta > 2$, corresponding to higher-order (greater than 2) moments among the $l'(\pi_{ij}^*)$ terms. Assumptions~IX(i) and IX(iii) are analogous to Assumption~VIII, but incorporate Hessian weights with $\zeta = 2$, whereas Assumption~IX(ii) resembles Assumption~VII(i) but concerns fourth-order moments among the $l'(\pi_{ij})$ terms. These assumptions hold when the edge variables are weakly  correlated. Below, we provide Remark~\ref{rem:example_dept} with examples of specific link functions and weak-correlation regimes under which the assumptions hold.

\begin{remark}\label{rem:example_dept} 
    To better understand the moments among the $l'(\pi_{ij}^*)$ terms, we consider two specific link functions as examples. Following Remark~\ref{rem:log1} in Section~\ref{sec:assum}, under a latent space model with the logistic link, we have $l'(\pi_{ij}^*) = A_{ij} - p_{ij}^*$ and $l''(\pi_{ij}^*) = p_{ij}^*(1 - p_{ij}^*)$. Thus, $\mathbb{E}[\,l'(\pi_{mi}^*) \, l'(\pi_{mj}^*)] = \mathrm{cov}(A_{mi}, A_{mj})$ corresponds to the covariance between edges $A_{mi}$ and $A_{mj}$, and different moments among the $l'(\pi_{ij}^*)$ terms correspond to different centered moments among the edge variables. Under linear low-rank decomposition-type models, we have $l'(\pi_{ij}^*) = (A_{ij} - \pi_{ij}^*)/\delta^2$ and $l''(\pi_{ij}^*) = -1/\delta^2$, where $A_{ij} \sim \mathcal{N}(\pi_{ij}^*, \delta^2)$. In this case, the moments of the $l'(\pi_{ij}^*)$ terms represent centered moments among the edge variables, scaled by their variance $\delta^2$.


We now introduce some sufficient conditions under which the additional dependencies are weak enough to satisfy these assumptions. Consider, for example, Assumption~VII(i), which bounds the scaled average second moments between the $l'(\pi_{ij})$ terms given the latent vectors. It suffices to ensure
$\sum_{m:m\ne i,j} \mathbb{E}\big(l'(\pi_{mi}^*)\,l'(\pi_{mj}^*)\big) \le \sqrt{n}\,M_1$ for some constant $M_1 > 0$ and for all pairs of $(i,j)\in[n]\times[n]$. This condition holds under several scenarios. First, if additional dependencies are scattered, meaning only $O(\sqrt{n})$ out of the $n-2$ ($n-1$ if $i=j$) nodes in $[n]\setminus\{i,j\}$ have non-zero $\mathbb{E}\big(l'(\pi_{mi}^*)\,l'(\pi_{mj}^*)\big)$, the assumption naturally holds. In the latent space, this can be represented by sparse embeddings in the latent dimensions that are not estimated. Second, if the moments uniformly shrink, i.e.,
$\mathbb{E}\big(l'(\pi_{mi}^*)\,l'(\pi_{mj}^*)\big) \le {M_2}/{\sqrt{n}}$ for all $i,j,m$ and some constant $M_2 > 0$, the assumption also holds. The decaying-product correlation structure in \cite{jiang2021set} provides another example that meets this requirement. Similar sufficient conditions can be introduced for other moments among the $l'(\pi_{ij}^*)$ terms. In practice, these assumptions can be checked using their plug-in moment estimators.
In the simulation studies of Section~\ref{sec:simulation}, we consider specific dependent data settings to study the corresponding correlation settings aforementioned. We further investigate how to check the strength of additional dependency in the assumptions and its relationship to inference validity in Section~A.2 of the Supplementary Material. 
\end{remark}
}
We summarize the main theorems as follows.
{
\begin{theorem}[Uniform and Average Consistency]
\label{thm:uni_consistency_dep}
Under Assumptions II-VIII, we have 
$$\|\hat \phi - \phi^*\|_{\max} = O_p(n^{-1/2+1/\zeta})~\text{ and }~n^{-1}\|\hat \phi - \phi^*\|^2 = O_p(n^{-1}).$$
\end{theorem}
}
\begin{remark}
    Compared to Theorem \ref{thm:uni_consistency_bounded}, the uniform consistency rate becomes $O_p(n^{-1/2+1/\zeta})$, which relies on $\zeta$ defined in Assumption VIII. This indicates that the rate now depends on both the specific form of the likelihood function and the level of edge dependency. On the other hand, The average consistency rate remains optimal and consistent with the result in Section \ref{sec:bounded_case}.
\end{remark}
{
To characterize the asymptotic distributions under the dependent-edge setting, we introduce the following consistent estimators of covariance matrices introduced in Assumption X:
\begin{equation}\label{eq:variance-estimator2}
        \hat{\Omega}_\mathcal{I} = \Big([\hat \Omega_\mathcal{I}]_{ij}\Big)~ \text{with} ~
      [\hat \Omega_\mathcal{I}]_{ij} = \frac{1}{n}\sum_{k_1: k_1\neq i}\sum_{k_2: k_2 \neq j} l'(\hat \pi_{{ik_{1}}})l'(\hat \pi_{{jk_{2}}})\hat h_{k_1}{\hat h_{k_2}}^{\top}, ~ i, j \in \mathcal{I}.
\end{equation}}
\begin{theorem}[Joint Asymptotic Distribution]
\label{thm:joint_distribution_dep}
Given $m$ node indices $\mathcal{I} = (i_1, i_2, ..., i_m)$, denote $\phi_{\mathcal{I}} = (\phi_{i_1}^{\top}, \phi_{i_2}^{\top},...,\phi_{i_m}^{\top})^{\top}$. Under Assumptions II-X, we have
$$\sqrt n \Big[\widehat{\text{Var}(\hat \phi_{\mathcal{I}})}\Big]^{-1/2} (\hat \phi_{\mathcal{I}} - \phi_{\mathcal{I}}^*) \overset{d}{\to} \mathcal{N} (0, I_{m(r+1)}),$$
where $\widehat{\text{Var}(\hat \phi_\mathcal{I})} = \hat\Sigma_\mathcal{I}^{-1}\hat\Omega_\mathcal{I}\hat\Sigma_\mathcal{I}^{-1}$ and $\hat \Sigma_\mathcal{I}$ and $\hat \Omega_\mathcal{I}$ are defined in Equations \eqref{eq:variance_estimators} and \eqref{eq:variance-estimator2}.
\end{theorem}

\begin{remark}
  {The asymptotic variances under the dependent-edge setup have a ``sandwich'' form due to the additional dependency. Under the independent-edge setup, these variance terms themselves reduce to those established in Section \ref{sec:bounded_case}, as we show in Section B.6 of the Supplementary Material.} Similarly to Corollaries \ref{thm:ind_distribution_bounded} and Equation \eqref{cor:link_prediction_bounded}, the individual asymptotic distribution of every latent position estimator as well as the link prediction inference results hold under Assumptions II-X. Compared to existing literature,  asymptotic distribution results have only been established under the RDPG model \citep{athreya2017statistical} with the independent edge setting. On the other hand, our approach is more general and can be adapted to the dependent edge setting. 
\end{remark}

\subsection{Results for Sparse Networks}
\label{sec:sparse_results}
In the previous analysis under Assumption II, we only consider the bounded parameter space for the true parameters. This limitation primarily arises from the difficulty in the theoretical analysis, where the objective function is shown to exhibit strict concavity locally near the true parameters. 
However, this assumption poses an inherent restriction on the sparsity or signal-to-noise ratio of network data,  which may be impractical in certain applications.
For example, many social networks and communication networks are observed to have sparsity properties, indicating that the edge density tends to decrease as $n$ increases.


In this section, we illustrate how to relax Assumption II and conduct the theoretical analysis of the likelihood function within an unbounded parameter space. 
Our approach is to generalize the model setup in Equation \eqref{eq:general_model} by introducing a sparsity parameter that controls the overall edge density and is allowed to be unbounded. This new parameter allows us to quantify the influence of network sparsity on the properties of the maximum likelihood estimators. 
By modifying the penalty terms $P(\phi)$ and changing the assumptions accordingly, we are able to utilize the same analysis techniques to derive meaningful inference results. 
However, it is crucial to note that these modifications are case-specific and may vary for different link functions and edge types.
Consequently, we focus our analysis on the popular latent space model for binary networks,
with the logistic link function,  $\sigma(x) = \text{logistic}(x) = 1/(1 + \exp(-x))$, and Bernoulli edge variables. The extensions towards other link functions and edge types can be performed similarly; see Remark \ref{Gaussian-sparse} for more discussion.
 


Specifically, we generalize the the model in Equation \eqref{eq:general_model} into  Equation \eqref{eq:sparse_model} below, where $\rho_n \in \mathbb R$ is the newly introduced sparsity parameter that is allowed to be unbounded while $z_i$ and $\alpha_i$ remain bounded.
\begin{equation}
    A_{ij} \sim \text{Bernoulli}(\theta_{ij}), \quad  \theta_{ij} = \text{logistic}(z_i^{\top}z_j + \alpha_i + \alpha_j + \rho_n),
    \label{eq:sparse_model}
\end{equation}
Our theoretical analysis is conducted under the setting that the sparsity parameter $\rho_n$ is given, similar to the setting considered in the sparse low-rank matrix completion literature \citep{cape2019signal}. The estimators of $Z$ and $\alpha$ are defined similarly as in Section \ref{sec:problem_setup_notations}, that is, the corresponding maximum likelihood estimators under the boundedness and identifiability constraints for $Z$ and $\alpha$.
We emphasize that introducing $\rho_n$ serves solely the purpose of characterizing the effect of network sparsity on inference results and do not change our inference targets ($Z$ and $\alpha$). Similar to the treatment in \cite{cape2019signal}, $\rho_n$ shall be seen as a known parameter to characterize the effect of sparsity, on which the estimators' consistency rates would depend. 




To make the notations consistent with previous sections, we keep using $\pi_{ij}$ to denote $z_i^{\top} z_j + \alpha_i + \alpha_j$ and $\phi_i$ to denote $(z_i^{\top}, \alpha_i)^{\top}$.
Note that all the parameters of interests, $\phi_i$ and $\pi_{ij}$, do not include $\rho_n$. 
By separating the sparsity effect from the formulation of $\theta_{ij}$ with a standalone sparsity parameter ($\rho_n$), it facilitates our theoretical analysis of the maximum likelihood estimators for $Z$ and $\alpha$ as follows.
In particular, we follow the theoretical framework introduced in Section \ref{sec:bounded_case}. Writing out the objective function in Equation \eqref{eq:loss_function_bounded} in the same way, our theoretical analysis will focus on the following loss function:
\begin{equation}
\begin{split}
    Q(\phi) = L(\phi) + P(\phi) = & \sum_{i = 1}^n \sum_{j = i + 1}^n l(z_i^{\top}z_j + \alpha_i + \alpha_j + \rho_n) \\
    &- \frac{c_n}{2}n^2 \left\|\text{upper\_diag}(Z^{\top}Z/n)\right\|_F^2 - \frac{c_n}{2}n^2\Big\|\frac{Z^{\top}1_n}{n}\Big\|_F^2,
    \label{eq:loss_function_sparse}
\end{split}
\end{equation}
where we write the multiplier $c$ as $c_n$, as it will be chosen to depend on the network sparsity ($\rho_n$), whose rationale is explained as follows.  For simplicity, we still denote the likelihood function as $l(\pi_{ij}) =  l(z_i^{\top}z_j + \alpha_i + \alpha_j + \rho_n)$ when there is no ambiguity. 
Specifically, to address the challenge of strict concavity of $l(\pi_{ij})$ in the unbounded parameter space, we examine the limiting relations of $l'(\pi_{ij})$, $l''(\pi_{ij})$, $l'''(\pi_{ij})$ and $\rho_n$, so that after proper scaling the analysis reduces back into a bounded parameter space. For the   model in Equation \eqref{eq:sparse_model}, we have: 
\begin{equation}
\begin{split}
    l'(\pi_{ij}) &= A_{ij} - \theta_{ij}=:A_{ij}-p_{ij}, \\
    l''(\pi_{ij}) &= -\frac{e^{\pi_{ij}+\rho_n}}{(1 + e^{\pi_{ij}+\rho_n})^2}, \\
    l'''_{ij}(\pi_{ij}) &= -\frac{e^{\pi_{ij}+\rho_n}(-1 + e^{\pi_{ij}+\rho_n})}{(1 + e^{\pi_{ij}+\rho_n})^3}.
    \label{eq:derivatives}
\end{split}
\end{equation}
When $\rho_n \to -\infty$ and $\pi_{ij}$ remains bounded, all $l'(\pi_{ij})$, $l''(\pi_{ij})$, $l'''(\pi_{ij})$ can be rescaled with {$w_n:= e^{\rho_n}$} such that similar conditions stated in Assumption V are valid. Correspondingly, the scaling factor for the Hessian matrix will also become $(w_nn)^{-1}$ and the order of penalty term $P(\phi)$ in Equation \eqref{eq:loss_function_sparse} should multiply $w_n$. Thus in \eqref{eq:loss_function_sparse}, we   choose {$c_n = cw_n$} with a constant $0 < c < e^{-M}$, where $M$ is specified in the following Assumption II*. The theoretical properties of the estimators are established in the following.



\begin{remark}
    In practice when $\rho_n$ is unknown, our estimation approach is to first estimate $Z$ and $\alpha$ using the same procedure as in the independent-edge case of Section 2 (with the identifiability constraints), and then estimate $\rho_n$ as the mean of the estimated $\alpha$ in the first step and set the final estimator of $\alpha$ to be centered. Note that the centering for $\alpha$ is needed for identifiability consideration since the mean of $\alpha_i$'s and $\rho_n$ cannot be distinguished when $\rho_n$ is unknown.  Empirical consistency results are presented in the simulation results under the ``Sparse'' setting in Section \ref{sec:simulation}.
    As will be shown below in our theoretical results, constructing confidence regions or intervals does not directly require the estimation of $\rho_n$, and our numerical results show that the established theory 
remains valid in the considered sparse edge setting.
\end{remark}




\begin{remark}\label{Gaussian-sparse}
    Before presenting the theoretical results, we briefly discuss how similar modifications can be made for other link functions and edge types. The main idea is to introduce a proper scaling parameter for the log-likelihood function of the original model, $l(\pi_{ij})$, and modify the assumptions based on the analysis of the order of its derivatives. Take the Gaussian edge variable with identity function as an example, it is reasonable to introduce a multiplicative parameter and consider the model 
\begin{equation*}
    A_{ij} \sim \mathcal N\left(\theta_{ij}, \delta^2\right), \quad  \theta_{ij} = \rho_n \sigma (z_i^{\top}z_j + \alpha_i + \alpha_j ),
\end{equation*}
where $\rho_n$, in this case, controls the overall signal-to-noise ratio and is allowed to decrease to zero for modeling low signal-to-noise ratio data.
\end{remark}

We make the following assumptions to analyze the model in Equation \eqref{eq:sparse_model} under the sparse and dependent-edge network setting.  

{
\begin{enumerate}[I*.]
    \item With $w_n = e^{\rho_n}$, we have $\lim_{n\to \infty} w_n = 0$. Furthermore, there exists $\epsilon > 0$ such that $w_n = \Theta(n^{-1/2 + \epsilon})$, i.e., there exists $C > 0$ such that $w_n \geq C n^{-1/2 + \epsilon}$.
\smallskip
    \item There exists $M > 0$ such that the parameter space $\Xi = \left\{\phi | \sup_{1 \leq i \leq n} \|\phi_i\| \leq M \right\}$ contains  the true parameters $\phi^*$. 
    \smallskip
    \item There exists a positive definite $r \times r$ matrix $\Sigma_Z$ such that ${Z^*}^{\top}Z^*/n \to \Sigma_Z$ as $n \to \infty$. 
\smallskip
    \item The limiting covariance matrix $\Sigma_Z$ is diagonal with with unique eigenvalues. Also, $Z^*$ is centered such that ${Z^*}^{\top}1_n = 0_r$.
 \smallskip   
    \item There exists $M > 0$ large enough such that
    \begin{enumerate}[i.]
        \item $n^{-2}\sum\limits_{i=1}^n \sum\limits_{j = 1}^n \big\{\frac{1}{w_n\sqrt{n}}\sum\limits_{m:m\ne i,j}\mathrm{cov}(A_{mi},A_{mj})       \big\}^2 \leq M.$
        \item {$n^{-2}\sum\limits_{i=1}^n \sum\limits_{j = 1}^n\mathbb{E}\big[\frac{1}{w_n\sqrt n}\sum\limits_{m:m\ne i,j}\big\{(A_{mi}-p_{mi}^*)(A_{mj}-p^*_{mj}) - \mathrm{cov}(A_{mi},A_{mj}) \big\}\big]^2 \leq M$.} 
    \end{enumerate}
    \item There exists $M > 0$ large enough and $\zeta > 2/(\epsilon+2^{-1})$ such that\\
        $\mathbb{E}({n}^{-1}\sum\limits_{i=1}^n\|\frac{1}{\sqrt{w_nn}}\sum\limits_{j: j \neq i} (A_{ij}-p_{ij}^*)h^*_j\|^{\zeta}) \leq M$, and  $\EE(\max\limits_{i\in[n]}\sum\limits_{j:j\ne i} \mathbf{1}_{\{A_{ij}=1\}} ) \le Mw_nn$.
    \item There exists $M > 0$ large enough such that for all $1 \leq s \leq n$ and $1\le q\le r+1$: 
    \smallskip
    \begin{enumerate}[i.]
        \item $\EE\big\|{n}^{-1}\sum\limits_{j:j\ne s}\sum\limits_{m:m\ne j}\big(n^{-1}\sum\limits_{l:l\ne j}l''(\pi_{jl}^*)h_l^* {h_l^*}^{\top}\big)^{-1}(A_{jm}-p_{jm}^*)h_m^* {h_j^*}^{\top}l''(\pi_{si}^*)\big\| \leq M,$
        \item {$\EE\big\|{n}^{-1}\sum\limits_{j:j\neq s}\sum\limits_{m:m\neq j}\big(n^{-1}\sum\limits_{l:l\ne j}l''(\pi_{jl}^*)h_l^*{h_l^*}^{\top}\big)^{-1}(A_{sj}-p_{sj}^*)(A_{jm}-p_{jm}^*)h_m^*\big\| \leq M$,}
        \item $\EE\big\|n^{-1}\sum\limits_{j=1}^n\sum\limits_{m: m \neq j}h_{jq}^*\big((nw_n)^{-1}\sum\limits_{l: l\neq j} l''(\pi_{jl}^*)h_l^*{h_l}^{*\top}\big)^{-1} w_n^{-1/2} (A_{jm} - p_{jm}^*) h_m^*\big\|\le M$.
    \end{enumerate}
    \item For an integer $m$ and any $m$ node indices $\mathcal{I} = (i_1, i_2, ..., i_m)$, there exists a non-random sequence  $\Omega_{\mathcal{I}}:=\Omega_{n,\cI}$ such that $(w_nn)^{-1/2}\Omega_{\cI}^{-1/2}\big(\sum_{j:j\ne i_1}(A_{i_1j}-p^*_{i_1j})h_j^*, \sum_{j:j\ne i_2}\\(A_{i_2j}-p^*_{i_2j})h_j^*, ..., \sum_{j:j\ne i_1}(A_{i_mj}-p^*_{i_mj})h_j^*\big) \overset{d}{\to} \mathcal{N}(0, I_{m(r+1)})$ as $n\to\infty$.
\end{enumerate}
}

Our analysis for the sparse network setting under \eqref{eq:sparse_model} is compatible with the dependent edge setting, where the orders of the quantities in Assumptions V*-VIII* are rescaled with $w_n$ accordingly. {In particular, the second part of Assumption VI* is an additional assumption concerning the strength of edge dependence, introduced to account for sparsity in the dependent-edge case. When the edges are independent, we verify this condition by bounding the maximum node degree in an Erdős--Rényi random graph \citep{erdos1960evolution}, as discussed in Section~B.6 of the Supplementary Material. When the edges are not independent, this assumption holds under weak edge dependence conditions, such as tail bounds on the node degrees \( \{\sum_{j: j \ne i} \mathbf{1}_{\{A_{ij} = 1\}} \}_{i=1}^n \) and structural assumptions on the edge dependencies \citep{janson2004large, chatterjee2007stein}.} Assumption I* requires the edge density of the network to be at least of order $n^{-1/2 + \epsilon}$, such that the effective sample size $w_n n$ is greater than $\sqrt{n}$ when allowing for dependency. {Assumptions II*-IV* are similar to the previous Assumptions II-IV on model regularity and identifiability.} The previous Assumptions V-VI on the likelihood function are satisfied by the model in Equation \eqref{eq:sparse_model}, hence they are not needed.

{
\begin{remark}\label{rem:spar_dep1} 
When proving uniform consistency and asymptotic normality of the maximum likelihood estimators, it is essential to understand the high-dimensional structures of the score vector and the Hessian matrix. These structures differ substantially between dependent-edge and independent-edge cases. For dependent-edge setups, classical random matrix concentration techniques that rely on independence among entries are not applicable \citep{chung2011spectra,lei2015consistency,paul2020spectral}. 
To characterize the behavior of the score and Hessian under dependent-edge and sparse regimes, we introduce Assumptions~V*--VII*, which impose bounds on the strength of additional dependency in a sparse setting. Rather than fixing a specific parametric form of the additional dependency, we employ general moment conditions, broadening the applicability of our theoretical results. These assumptions not only serve technical purposes in our theoretical derivations, but also provide practical diagnostics for assessing the validity of inference, as illustrated in Section~A.2 of the Supplementary Material. These conditions also offer insights into the trade-off between sparsity and dependency: the sparser the network, the weaker the additional dependency that can be tolerated. 
Under such dependent-edge setups, it is difficult to determine the optimality of the sparsity requirement in Assumption~I*, since we do not specify the exact form of the additional dependency but instead make general assumptions about its magnitude. When Assumption~I holds, i.e., when all edge variables are conditionally independent, we utilize different strategies to analyze the high-dimensional properties of the score and Hessian, allowing for a relaxation of the sparsity condition. In particular, we demonstrate later in this section that Assumption~I* can be weakened to permit network sparsity as low as $\Theta(n^{-1+\epsilon})$ under independent-edge setups. This minimal sparsity level aligns with the existing literature on inference for RDPG \citep{athreya2017statistical}, which is a specific type of latent space models with an identity link function for binary independent edges, where an arbitrarily slow polynomial term $n^{\epsilon}$ replaces the $(\log n)^{4+\epsilon}$ term therein. In contrast, our analysis covers general classes of latent space models (both linear and non-linear), providing analytic tools for general network data across a broad range of applications.
\end{remark}
}

Our theoretical results are summarized as follows.

{
\begin{theorem}[Uniform and Average Consistency]
\label{thm:uni_consistency_sparse}
Under Assumptions I*-VI*, we have 
$$\|\hat \phi - \phi^*\|_{\max} = O_p( w_n^{-1/2} n^{-1/2+1/\zeta}  )~~\text{and}~~\frac{1}{n}\|\hat \phi - \phi^*\|^2 = O_p\Big(\frac{1}{w_nn}\Big).$$
\end{theorem}
}

\begin{remark}
    In the sparse network setting, the consistency rates of maximum likelihood estimators are related to the rate at which $w_n$ approaches $0$. Similar results have been reported in the literature on sparse low-rank matrix completion \citep{cape2019signal} with a signal-plus-noise linear structure. However, our analysis adopts different techniques that are general and compatible with the dependent edge setting and non-linear model setups.
\end{remark}

Before presenting the results on statistical inference under the sparse edge setting, note that the consistent estimators of covariance matrices are similar but not the same as those presented in Equations \eqref{eq:variance_estimators} and \eqref{eq:variance-estimator2}. Given $m$ node indices $\mathcal{I} = (i_1, i_2, ..., i_m)$, the following consistent estimators are scaled with the sparsity parameter $w_n$: 

\begin{equation}
    \label{eq:variance_estimators_sparse}
    \begin{split}
        \hat \Sigma_{\mathcal{I}} &= \text{diag}(\hat \Sigma_{i_1}, \hat\Sigma_{i_2},...,\hat\Sigma_{i_m}), ~ \text{with} ~
        \hat \Sigma_{i} = - \frac{1}{w_nn}\sum_{j: j \neq i} l''(\hat \pi_{ij}) \hat h_j \hat h_j^{\top}, ~ i \in \mathcal{I}; \\ 
        \hat \Omega_\mathcal{I} &= \Big([\hat \Omega_\mathcal{I}]_{ij}\Big),~ \text{with} ~
        [\hat \Omega_\mathcal{I}]_{ij} = \frac{1}{w_nn}\sum_{k_1: k_1\neq i}\sum_{k_2: k_2 \neq j} l'(\hat \pi_{{ik_{1}}})l'(\hat \pi_{{ik_{2}}})\hat h_{k_1}{\hat h_{k_2}}^{\top}, ~ i, j \in \mathcal{I}.
    \end{split}
\end{equation}

\begin{theorem}[Joint Asymptotic Distribution]
\label{thm:joint_distribution_sparse}
Given $m$ node indices $\mathcal{I} = (i_1, i_2, ..., i_m)$, denote $\phi_{\mathcal{I}} = (\phi_{i_1}^{\top}, \phi_{i_2}^{\top},...,\phi_{i_m}^{\top})^{\top}$. Under Assumptions I*-VIII*,   we have
$$\sqrt{w_nn} \Big[\widehat{\text{Var}(\hat \phi_{\mathcal{I}})}\Big]^{-1/2} (\hat \phi_{\mathcal{I}} - \phi_{\mathcal{I}}^*) \overset{d}{\to} \mathcal{N} (0, I_{m(r+1)}),$$
where $\widehat{\text{Var}(\hat \phi_\mathcal{I})} = \hat\Sigma_\mathcal{I}^{-1}\hat\Omega_\mathcal{I}\hat\Sigma_\mathcal{I}^{-1}$ and $\hat \Sigma_\mathcal{I}$ together with $\hat \Omega_\mathcal{I}$ are defined in  \eqref{eq:variance_estimators_sparse}.
\end{theorem}

\begin{corollary}[Individual Asymptotic Distribution]
\label{thm:ind_distribution_sparse}
Under Assumptions I*-VIII*,   for each $1 \leq i \leq n$,  
$$\sqrt{w_nn} \Big[\widehat{\text{Var}(\hat \phi_i)}\Big]^{-1/2} (\hat \phi_i - \phi_i^*) \overset{d}{\to} \mathcal{N} (0, I_{r+1}),$$
where $\widehat{\text{Var}(\hat \phi_i)} = \hat\Sigma_i^{-1}\hat\Omega_i\hat\Sigma_i^{-1}$, with $\hat \Sigma_i$ and $\Omega_i$  defined in  \eqref{eq:variance_estimators_sparse}. 
\end{corollary}

    Importantly, Theorem \ref{thm:joint_distribution_sparse} and Corollary \ref{thm:ind_distribution_sparse} imply that constructing confidence regions and confidence intervals based on $\hat \phi_\mathcal{I}$  and $\hat \phi_i$, respectively,  does not require the estimation of $\rho_n$ or equivalently $w_n$. This is because (1) the scaling factor of the asymptotic variance term and the normalization factor in limiting distribution are canceled out; and (2) the asymptotic terms related to $w_n$ or $\rho_n$ in $l'(\pi_{ij})$ and $l''(\pi_{ij})$ are asymptotically canceled out.

Similar to Equation \eqref{cor:link_prediction_bounded}, we use the following example as an illustration for link prediction inference. Given $\mathcal{I} = \{i, j\}$, $1 \leq i < j \leq n$, $A_{ij}$ are Bernoulli random variables and $\sigma(x) = 1 / (1 + e^{-x})$, under Assumptions I*-VIII*, we have
\begin{equation}
    \label{cor:link_prediction_sparse}
    \sqrt{w_nn} \Big[\widehat{\text{Var}(\hat \theta_{ij})}\Big]^{-1/2} (\hat \theta_{ij} - \theta_{ij}^*) \overset{d}{\to} \mathcal N \left(0, 1\right),
\end{equation}
where in this case $\theta_{ij}$ is the linkage probability between node $i$ and node $j$, and
$$\widehat{\text{Var}(\hat \theta_{ij})} = \left[\hat\theta_{ij}(1-\hat\theta_{ij})\right]^2\begin{pmatrix} \hat h_j \\ \hat h_i  \end{pmatrix}^{\top} \begin{pmatrix}\hat\Sigma_i^{-1}  & 0 \\0 & \hat\Sigma^{-1}_j\end{pmatrix} \hat \Omega_{\mathcal{I}} \begin{pmatrix}\hat\Sigma_i^{-1}  & 0 \\0 & \hat\Sigma^{-1}_j\end{pmatrix}\begin{pmatrix} \hat h_j \\ \hat h_i \end{pmatrix}.$$
  Similar to  Theorem \ref{thm:joint_distribution_sparse} and Corollary \ref{thm:ind_distribution_sparse}, 
    constructing confidence regions or intervals for $\theta_{ij}$ does not require estimating $\rho_n$. When constructing intervals in practice, one could further remove the $w_n$ implicitly contained in $\theta^*_{ij}$ in the variance term with asymptotic approximations. However, through simulation studies, we find that constructing confidence intervals with the current form gives better coverage rates. This is potentially because removing $w_n$ would cause finite-sample estimation bias.

{
The previous results in this section concern sparse and dependent-edge setups. As Remark~\ref{rem:spar_dep1} indicates, considering independent-edge setups can further expand the sparsity regime in Assumption~I*. We now introduce the following assumption on model \eqref{eq:sparse_model} with an expanded sparsity regime, which leads to Theorem~\ref{thm:ind_spaese}.

\begin{enumerate}[I$^\dagger$.]
    \item With $w_n = e^{\rho_n}$, we have $\lim_{n\to \infty} w_n = 0$. Furthermore, there exists arbitrarily small  $\epsilon > 0$ such that $w_n = \Theta(n^{-1+\epsilon})$, i.e., there exists $C > 0$ such that $w_n \geq C n^{-1+\epsilon}$.
\end{enumerate}

\begin{theorem}\label{thm:ind_spaese}
    Under Assumptions I, I$^{\dagger}$, and II*-IV*, 
    we have 
    \[
      \|\hat{\phi} - \phi^*  \|_{\max} = O_p( w_n^{-1/2} n^{-1/2+\epsilon} )~\text{ and } ~n^{-1}\|\hat{\phi} - \phi^*\|^2 = O_p(w_n^{-1}n^{-1})
    \]
    for any $\varepsilon>0$. Moreover, given $m$ node indices $\mathcal{I} = (i_1, i_2, ..., i_m)$, let $\phi_{\mathcal{I}} = (\phi_{i_1}^{\top}, \phi_{i_2}^{\top},...,\phi_{i_m}^{\top})^{\top}$.   We have
$$\sqrt{w_nn} \Big[\widehat{\text{Var}(\hat \phi_{\mathcal{I}})}\Big]^{-1/2} (\hat \phi_{\mathcal{I}} - \phi_{\mathcal{I}}^*) \overset{d}{\to} \mathcal{N} (0, I_{m(r+1)}),$$
where $\widehat{\text{Var}(\hat \phi_\mathcal{I})} = \hat\Sigma_\mathcal{I}^{-1}$ and $\hat \Sigma_\mathcal{I}$ is defined in  \eqref{eq:variance_estimators_sparse}.
\end{theorem}

\begin{remark}\label{rem:sparse_res}
    Under the latent space model in Equation \eqref{eq:sparse_model} and Assumption~I, in which edges are conditionally independent given the latent vectors, Assumptions~V*--VIII* hold for any finite $\zeta>2$ in Assumption VI*, and Assumption~I* can be relaxed to Assumption~I$^\dagger$. The sparsity requirement in Assumption~I$^\dagger$ aligns with the existing literature~\citep{lei2015consistency,athreya2017statistical} up to an arbitrarily slow polynomial term, showing that our method adapts to the optimal sparsity regime under independent-edge setups.
Additionally, we are able to draw a similar conclusion as in Remark \ref{remark:spectra} when the eigenvalue uniqueness condition is not met, except that the bound \(O_p(n^{-{1}/{2}} \log^{{1}/{2}}n)\) is replaced by \(O_p(w_n^{-{1}/{2}}n^{-{1}/{2}} \log^{{1}/{2}}n)\) due to the sparsity of the network.
\end{remark}
}

\section{Simulation Studies}
\label{sec:simulation}

\subsection{Experiment Setup}

In this section, we present numerical experiment results on (1) evaluating consistency rates of the maximum likelihood estimators; (2) constructing confidence intervals for latent positions; and (3) constructing confidence intervals for link prediction probabilities. Here we focus on the latent space model for binary networks with link function $\sigma(x) = \text{logistic}(x) = 1/(1 + \exp(-x))$ and Bernoulli edge variables, as presented in Equation \eqref{eq:sparse_model}, under different simulation settings including independent edges with bounded parameter space, dependent edges, and sparse networks. 
Simulation results of similar experiments for the network latent space model with continuous Gaussian edge variables can be found in the Supplementary Material.
Details of the data generating processes will be discussed in each setting below. For the estimation procedure, we adopt a  singular value thresholding approach as in Algorithm 3 of \cite{ma2020universal}  and a projected gradient descent method as in Algorithm 1 of \cite{ma2020universal} for initialization and optimization, respectively.

\subsection{Independent Edges with Bounded Parameter Spaces}
\label{sec:simulation_ind_bounded}
We first consider the independent setting as discussed in Section \ref{sec:bounded_case}, referred to as ``Bounded \& Indep.''. The network data is generated as follows. We set $r = 2$, generate i.i.d. entries of $Z$ from a truncated normal distribution $\mathcal{N}_{[-2, 2]}(0, 1)$ and i.i.d. entries of $\alpha$ from a uniform distribution $U([1, 3])$. To satisfy the identifiability conditions, we first center $Z$ and $\alpha$, and set $\rho_n = -3$. Under such data generating process, the average edge density of networks is around $0.08$. The number of nodes $n$ ranges from 500 to 8000, and for each $n$, we repeat the experiment 200 times by sampling i.i.d. networks from the same $Z$, $\alpha$, and $\rho_n$.

\begin{figure}
    \centering
    \includegraphics[width = \linewidth]{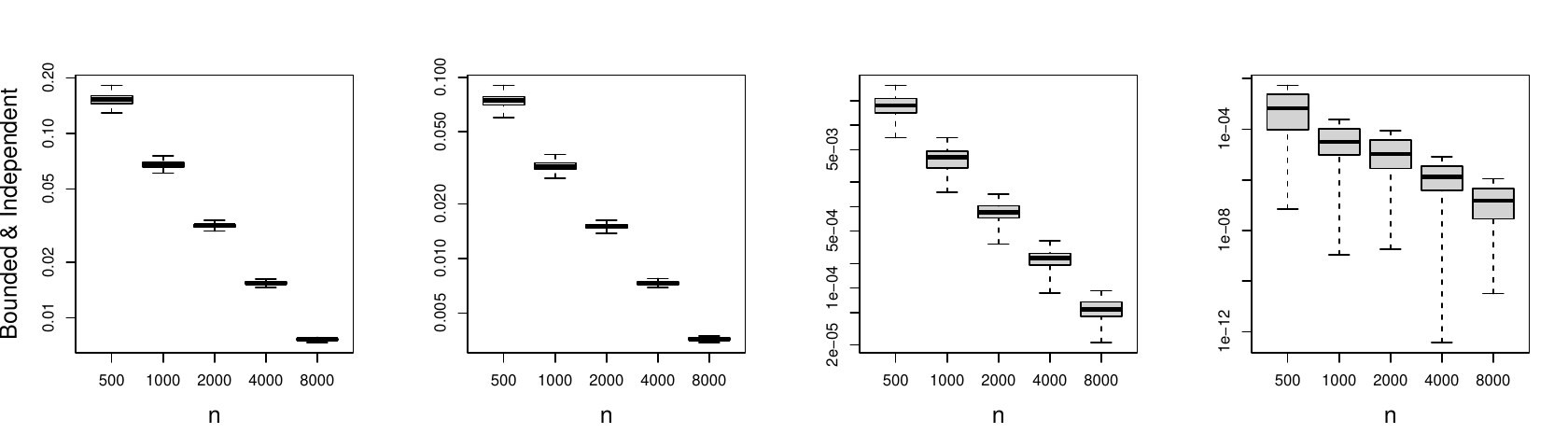}
    \includegraphics[width = \linewidth]{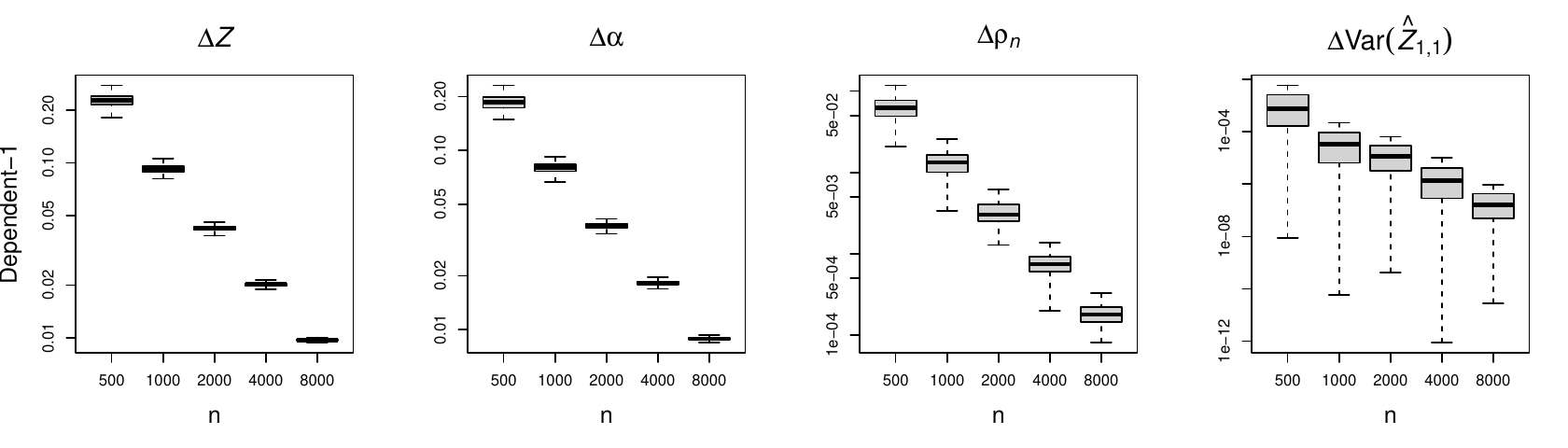}
    \includegraphics[width = \linewidth]{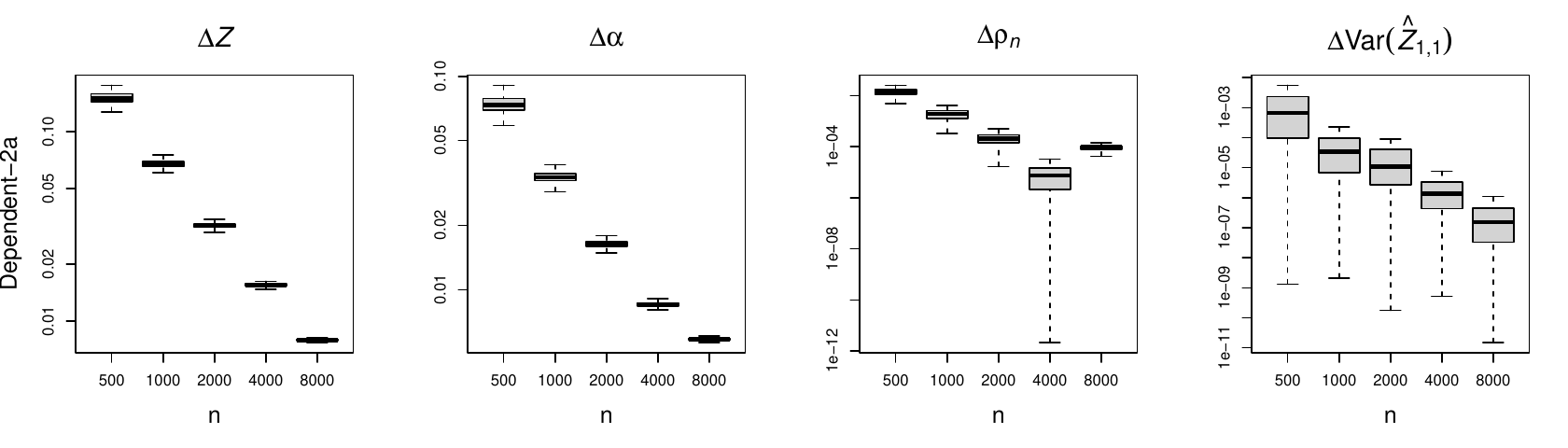}
    \includegraphics[width = \linewidth]{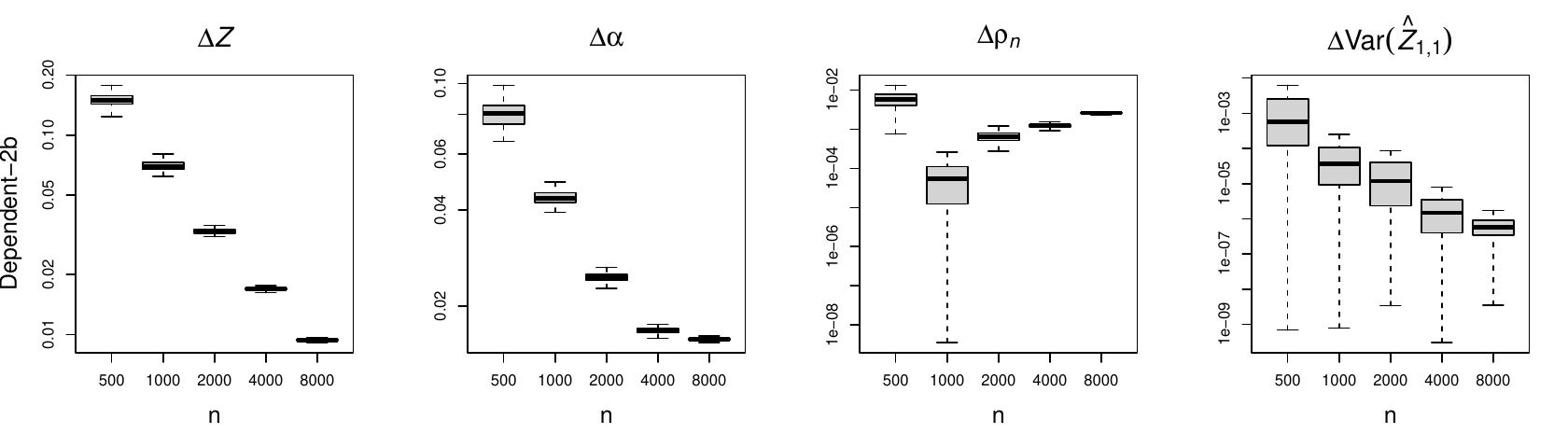}
    \includegraphics[width = \linewidth]{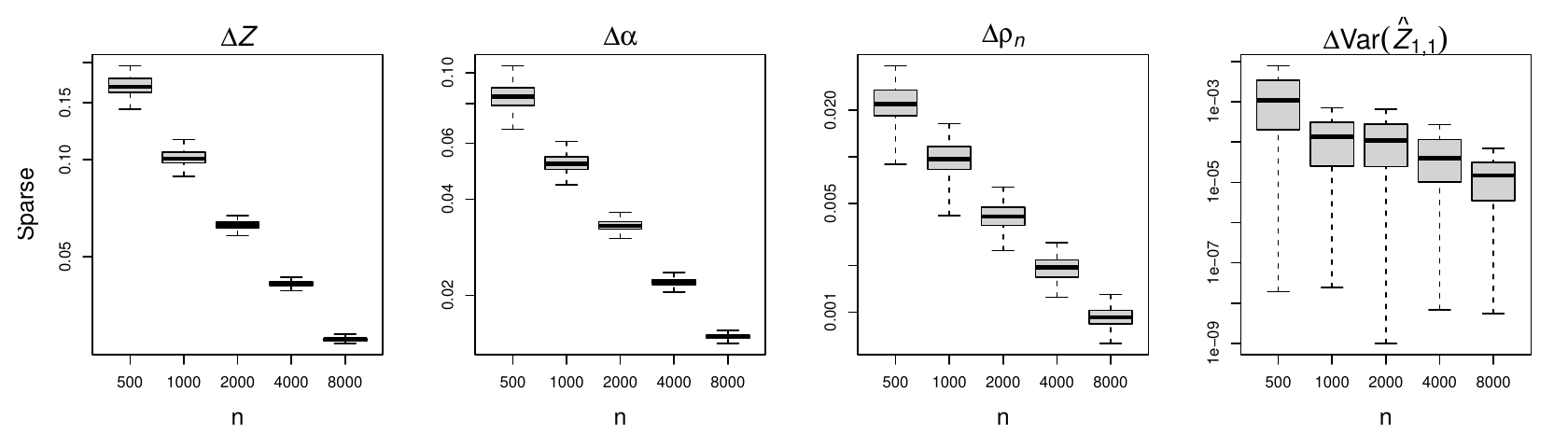}
    \caption{Consistency plots for maximum likelihood estimators under different settings.}
    \label{fig:consistency}
\end{figure}

For the consistency of maximum likelihood estimators, we evaluate $\Delta Z = \|Z - \hat Z\|_F^2 / n^2$, $\Delta \alpha = \|\alpha - \hat \alpha\|^2 / n^2$, $\Delta \rho_n = (\rho_n - \hat\rho_n)^2$, and $\Delta \text{Var}(\hat z_{11}) = (\text{Var}(\hat z_{11}) - \widehat{\text{Var}(\hat z_{11})})^2$, where the last quantity is measuring the estimation error of the plug-in estimator of the first coordinate of the first latent position (Equation \eqref{eq:variance_estimators}) and we estimate $\rho_n$ as the mean of $\hat \alpha$. 
The results are illustrated in the log-log scaled plots on the first row of Figure \ref{fig:consistency}. We can see the maximum likelihood estimators have consistency rates that meet theoretical results and both $\rho_n$ and $\text{Var}(\hat z_{11})$ can be consistently estimated.

\begin{figure}
    \centering
    \includegraphics[width = 0.95\linewidth]{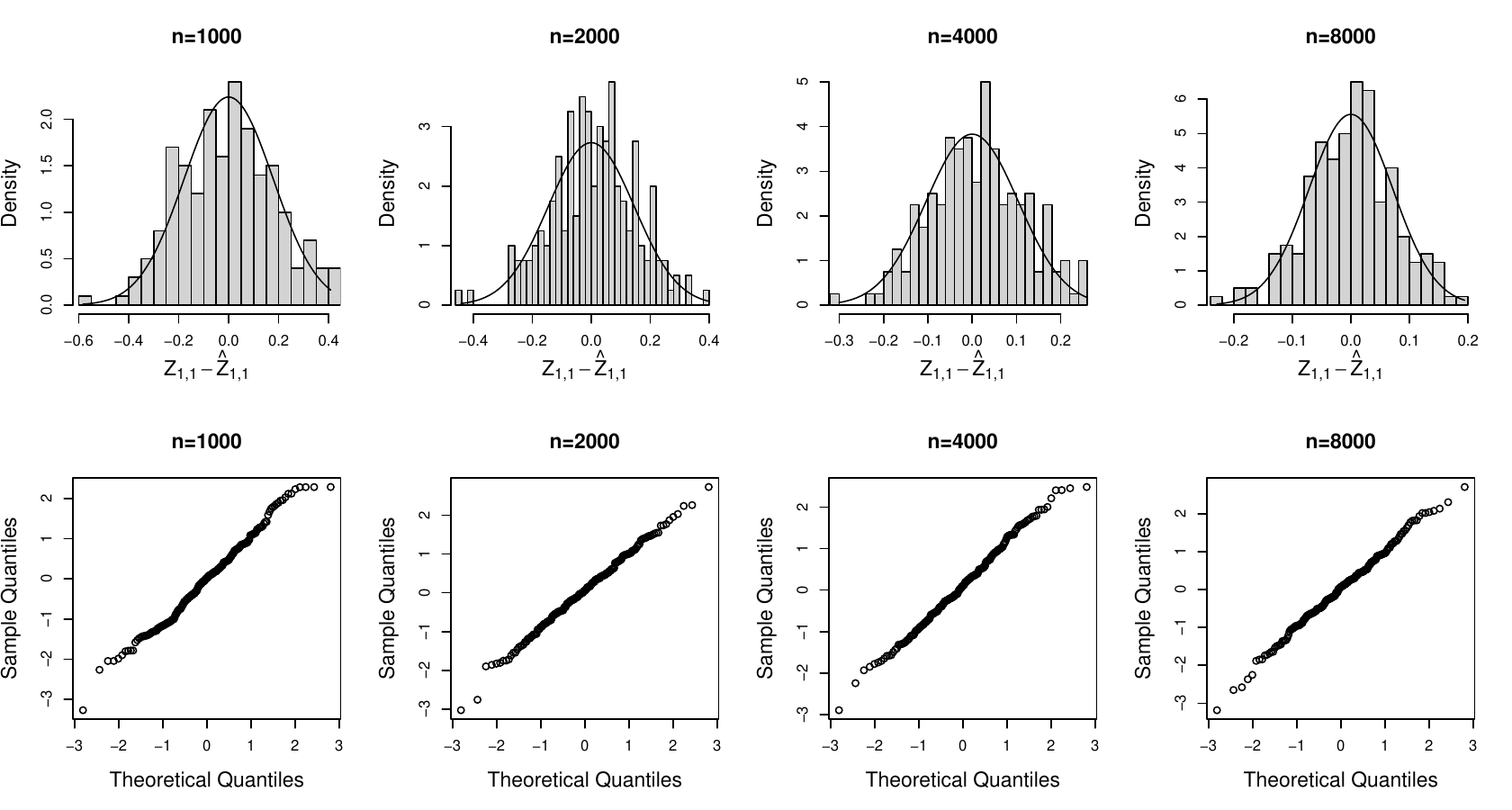}
    \caption{Empirical versus theoretical distribution plots under the ``Bounded \& Indep.'' setting.}
    \label{fig:distribution_plots_bounded}
\end{figure}

To evaluate the distributional results of latent positions, we focus on the first coordinate of the first latent position $z_{11}$ for an illustration. Figure \ref{fig:distribution_plots_bounded} plots the histograms of $z_{11} - \hat z_{11}$ together with the theoretical density curves, as well as the QQ-plots, where we can see that as $n$ grows the empirical distribution is converging to the theoretical distribution. We further construct 95\% confidence intervals for $z_{11}$ and the linkage probability of nodes 1 and 2, i.e., $\theta_{1,2}$, according to Corollary \ref{thm:ind_distribution_bounded} and Equation \eqref{cor:link_prediction_bounded} with the plug-in variance estimators. The average coverage rates are reported in the first rows of Table \ref{table:coverage_z11_w_std} and Table \ref{table:coverage_theta12_w_std}, respectively, where we can see that the empirical coverage rates meet the expectation as $n$ grows.

\begin{figure}
    \centering
     {\footnotesize Histograms and QQ-plots under Dependent-1}  \includegraphics[width = 0.95\linewidth]{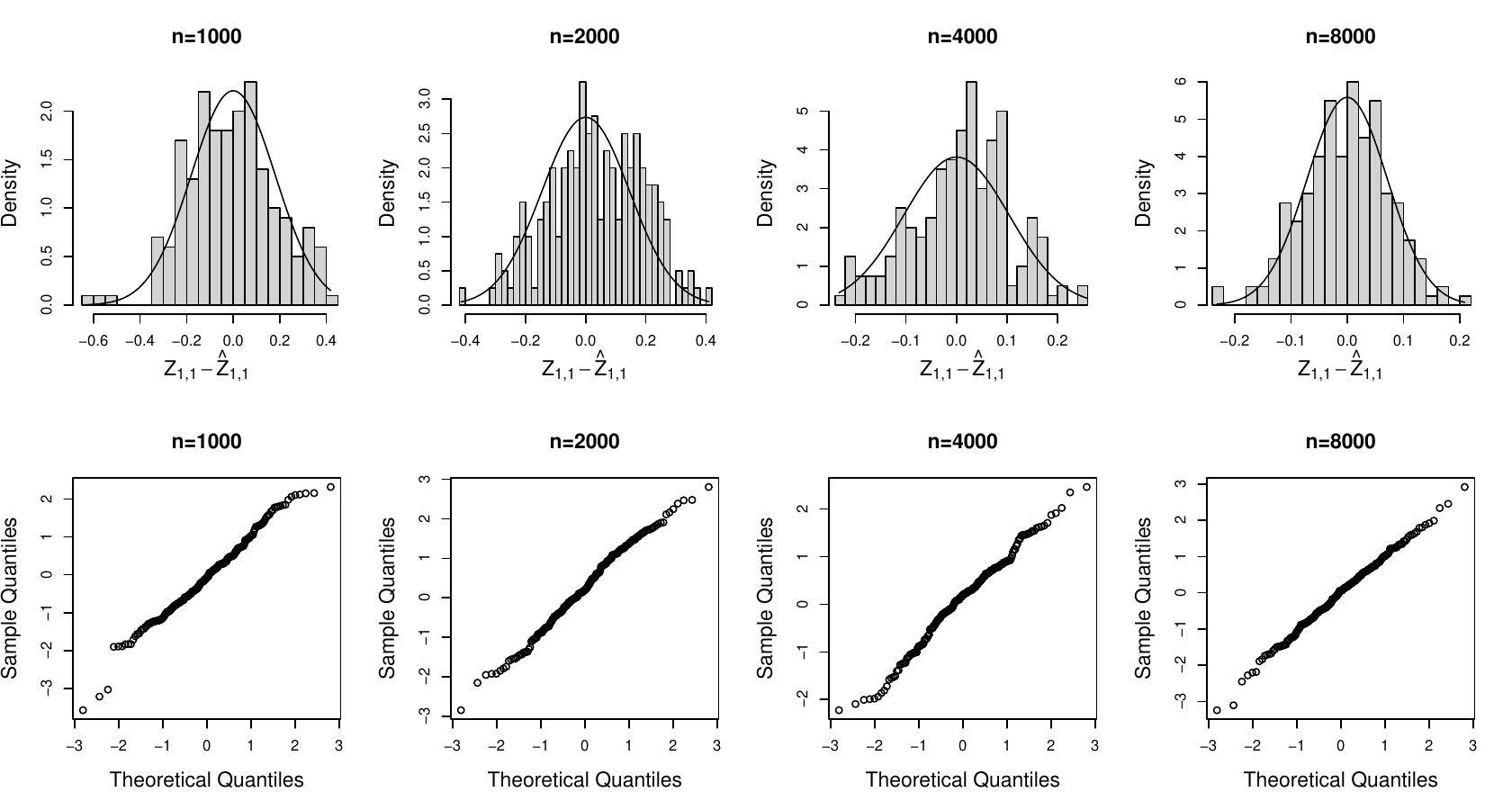}

   {\footnotesize  Histograms and QQ-plots under Dependent-2a}
    \includegraphics[width = 0.95\linewidth]{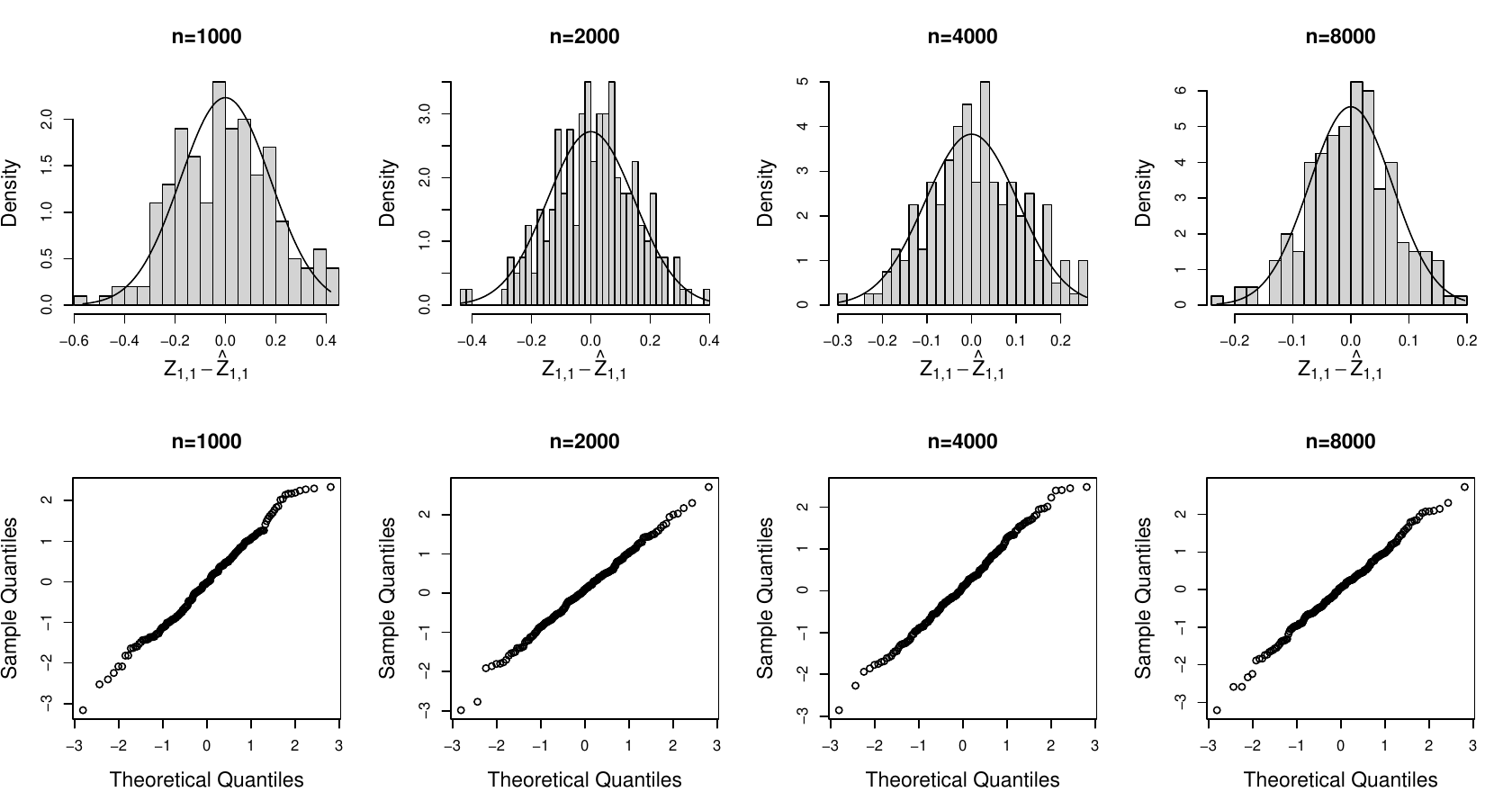}

       {\footnotesize  Histograms and QQ-plots under Dependent-2b}
 \includegraphics[width = 0.95\linewidth]{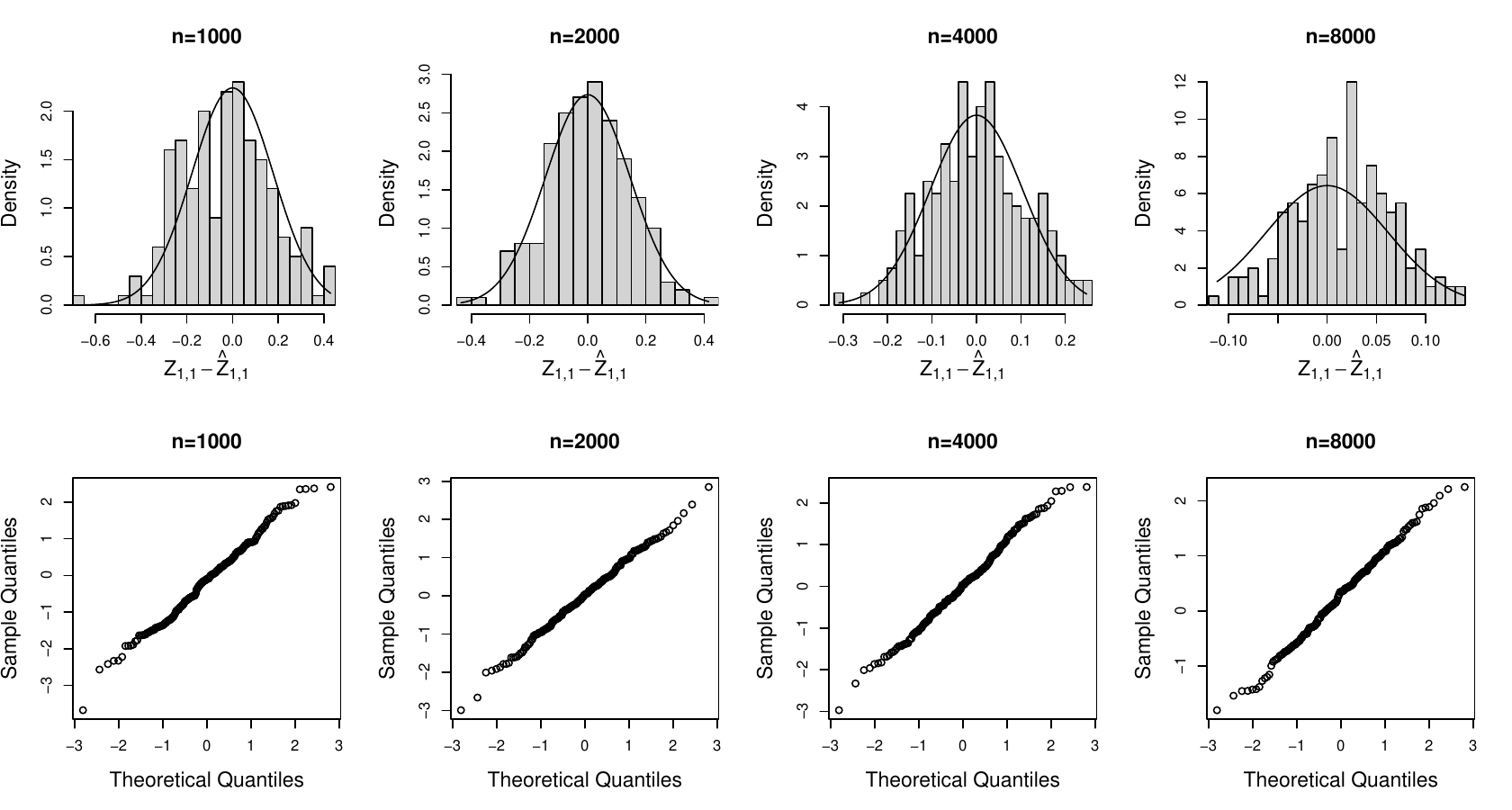}

    \caption{Empirical versus theoretical distribution plots under the dependent settings: ``Dependent-1'', ``Dependent-2a'', and ``Dependent-2b''.}
    \label{fig:distribution_plots_dependent_1}
\end{figure}


\begin{table}[ht]
\centering
\begin{tabular}{|c|ccccc|}
\hline
\textbf{Setting} & \multicolumn{5}{c|}{\textbf{Coverage of $z_{11}$}} \\
\hline
Bounded \& Indep. & 0.970 (0.0121) & 0.940 (0.0168) & 0.960 (0.0139) & 0.955 (0.0147) & 0.930 (0.0180) \\
Dependent-1 & 0.970 (0.0121) & 0.915 (0.0197) & 0.935 (0.0174) & 0.945 (0.0161) & 0.955 (0.0147) \\
Dependent-2a & 0.975 (0.0110) & 0.930 (0.0180) & 0.960 (0.0139) & 0.955 (0.0147) & 0.935 (0.0174) \\
Dependent-2b & 0.970 (0.0121) & 0.940 (0.0168) & 0.965 (0.0130) & 0.955 (0.0147) & 0.965 (0.0130) \\
Sparse & 0.975 (0.0110) & 0.930 (0.0180) & 0.940 (0.0168) & 0.945 (0.0161) & 0.930 (0.0180) \\
\hline
$n$ & 500 & 1000 & 2000 & 4000 & 8000 \\
\hline
\end{tabular}
\caption{Coverage rates of 95\% confidence intervals for $z_{11}$ with standard deviation under different settings.}
\label{table:coverage_z11_w_std}
\end{table}

\begin{table}[ht]
\centering
\begin{tabular}{|c|ccccc|}
\hline
\textbf{Setting} & \multicolumn{5}{c|}{\textbf{Coverage of $\theta_{12}$}} \\
\hline
Bounded \& Indep. & 0.890 (0.0221) & 0.940 (0.0168) & 0.955 (0.0147) & 0.925 (0.0186) & 0.935 (0.0174) \\
Dependent-1 & 0.905 (0.0207) & 0.925 (0.0186) & 0.875 (0.0234) & 0.925 (0.0186) & 0.940 (0.0168) \\
Dependent-2a & 0.885 (0.0226) & 0.940 (0.0168) & 0.950 (0.0154) & 0.925 (0.0186) & 0.940 (0.0168) \\
Dependent-2b & 0.885 (0.0226) & 0.960 (0.0139) & 0.930 (0.0180) & 0.865 (0.0242) & 0.000 (0.0000) \\
Sparse & 0.875 (0.0234) & 0.930 (0.0180) & 0.940 (0.0168) & 0.935 (0.0174) & 0.960 (0.0139) \\
\hline
$n$ & 500 & 1000 & 2000 & 4000 & 8000 \\
\hline
\end{tabular}
\caption{Coverage rates of 95\% confidence intervals for $\theta_{12}$ with standard deviation under different settings.}
\label{table:coverage_theta12_w_std}
\end{table}

    


\subsection{Dependent Edges}
\label{sec:dept-edge}
We consider two ways of introducing dependency between edges. In the first setting, referred as ``Dependent-1'', we generate $Z$, $\alpha$, and $\rho_n$ with the same setting as Section \ref{sec:simulation_ind_bounded}. Then we use R package ``CorBin" \citep{jiang2021set} to generate dependent Bernoulli edge variables with the decaying-product correlation structure \citep{jiang2021set}, where the correlation between nodes $i$ and $j$ are defined by $\prod _{l=i-1}^{j}\rho_{l}$ where $\rho_l$'s are determined by the marginal probabilities $\theta_{ij}$. 
The consistency plots are illustrated in the second row of Figure \ref{fig:consistency} and the coverage results are reported in the second rows of Table \ref{table:coverage_z11_w_std} and Table \ref{table:coverage_theta12_w_std}. 
The histograms and QQ-plots are illustrated in the first two rows of Figure \ref{fig:distribution_plots_dependent_1}. 
We can see that the theoretical results remain valid under this dependent setting.


The second setting being considered is to introduce dependency with an additional dimension of latent position that occurs in the data generating process while being ignored in model fitting. In addition to the generating process of $Z$, $\alpha$, and $\rho_n$ as before, we further generate 1-dimensional positions with Hardamard product $Z_{dep} \odot E$, where entries $Z_{dep}\overset{i.i.d.}{\sim} \mathcal{N}_{[-2,2]}(0, 1)$ and $E$ is an indicator vector with a certain proportion of random entries to be 1 and the rest being 0. The $Z_{dep}$ will only appear in data generation while not being estimated during model fitting, thus introducing implicit edge dependency. The non-zero proportions are set to be 25\% and 50\%, referred as ``Dependent-2a'' and ``Dependent-2b'' settings, respectively. Note that here we deliberately select relatively strong dependent cases with non-zero proportions 25\% and 50\% to give a full illustration of the strengths and the limitations of our theoretical results. 

The consistency plots are illustrated in the third and fourth rows of Figure \ref{fig:consistency},  the average coverage rates are reported in   Table \ref{table:coverage_z11_w_std} and Table \ref{table:coverage_theta12_w_std},  and the corresponding histograms and QQ-plots for   ``Dependent-2a'' and ``Dependent-2b'' settings are presented in Figure \ref{fig:distribution_plots_dependent_1}. 
From Figure \ref{fig:distribution_plots_dependent_1} and  Table \ref{table:coverage_z11_w_std}, we can see that the individual asymptotic approximation still performs well. On the other hand, it is interesting to observe from Figure \ref{fig:consistency} that $\hat \alpha$ and $\hat \rho_n$ are no longer showing ideal consistency patterns. Note that the scale of $\Delta \rho_n$ is much smaller than $\Delta \alpha$, and our theory focuses only on $\alpha$. Thus, the most serious deviation of the experiment result happens when $n = 8000$ under the ``Dependent-2b'' setting. To explain this observation, we check the validity of Assumption VIII and summarize the average values of $\|S(\hat \phi)\|$ of ``Dependent-1'', ``Dependent-2a'', and ``Dependent-2b'' with different values of $n$ in the Supplementary Material. Under Assumption VIII, the norm of the score vector should not increase significantly as $n$ increases. However, we find that $\|S(\hat \phi)\| = 0.53$ when $n = 8000$ under the ``Dependent-2b'' setting, whereas the values are around 0.13 for the rest situations. This indicates under the ``Dependent-2b'' setting, the dependency of edges might be too strong, and as a consequence the assumptions required for the validity of the theorems no longer hold. Interestingly, the coverage rates reported in the third and fourth rows of Table \ref{table:coverage_z11_w_std} and Table \ref{table:coverage_theta12_w_std} show consistent patterns. While most of the entries are around 0.95 as $n$ grows, when $n=8000$ under the ``Dependent-2b'' setting, the estimated variance of $\hat \theta_{1,2}$ is significantly smaller than the theoretical value due to the strong dependency, resulting in zero coverage. {We conduct additional numerical studies on dependent-edge settings and present the results in Section~A.2 of the Supplementary Material. }

\begin{figure}[ht]
    \centering
    \includegraphics[width = 0.95\linewidth]{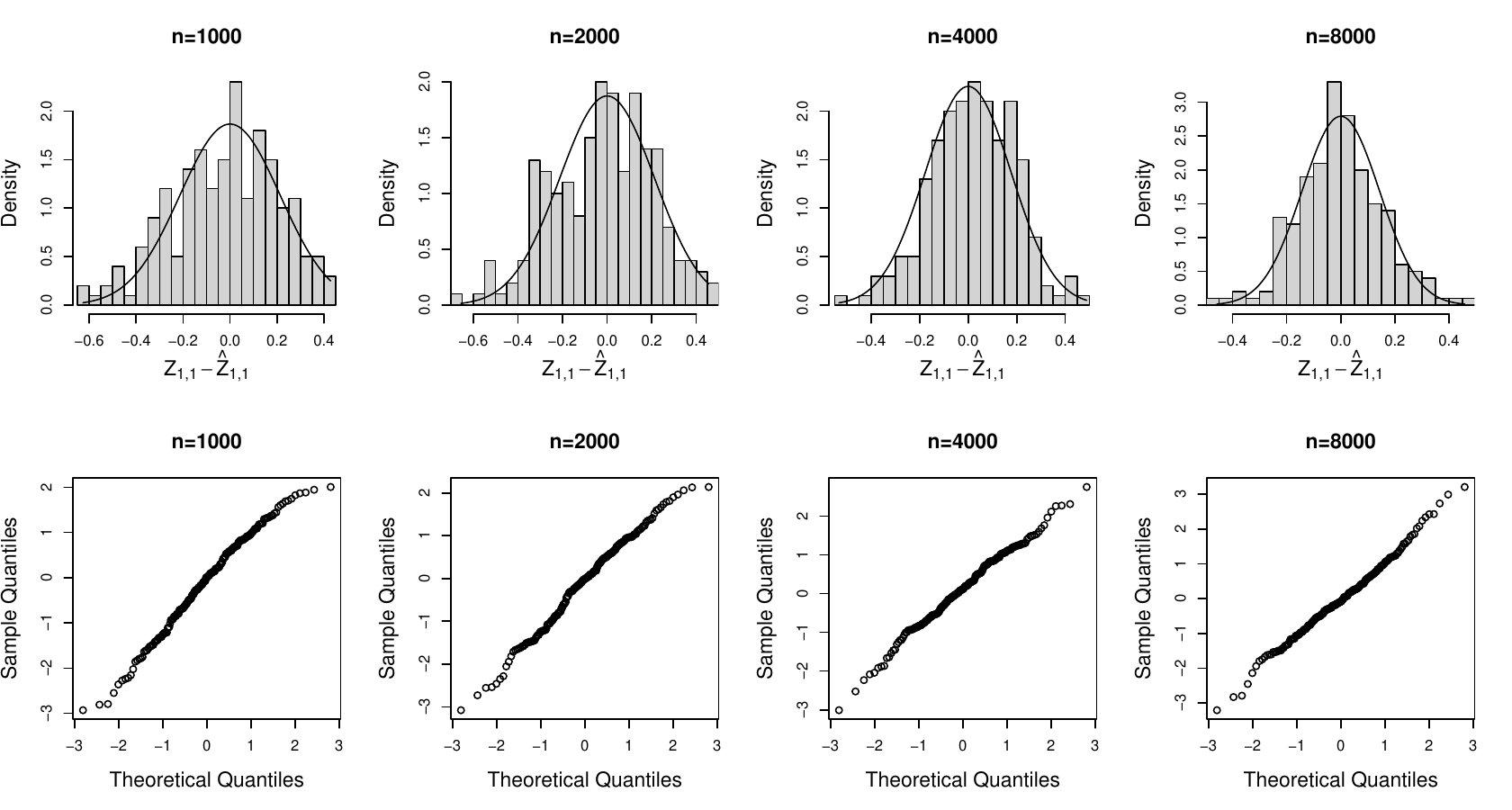}
    \caption{Empirical versus theoretical distribution plots under the ``Sparse'' setting.}
    \label{fig:distribution_plots_sparse}
\end{figure}

\subsection{Sparse Networks}

In the sparse edge setting, referred to as ``Sparse'', we generate $Z$ and $\alpha$ following the same way as in Section \ref{sec:simulation_ind_bounded} while setting $\rho_n = -\frac{1}{2}\log n$, which is on the borderline of violating Assumption I*. Under this construction, the average edge density of generated networks decreases from $0.07$ to $0.02$ as $n$ increases from $500$ to $8000$. The consistency plots are illustrated in the last row of Figure \ref{fig:consistency}, where we can see that the consistency rates is related to $\rho_n$ and meet the theoretical value of $O_p(\frac{1}{w_nn}) = O_p(\frac{1}{\sqrt{n}})$. The coverage rates are reported in the fifth rows of Table \ref{table:coverage_z11_w_std} and Table \ref{table:coverage_theta12_w_std}, showing that the empirical values meet the expectation well as $n$ grows. The histograms and QQ-plots for the distribution of latent position estimators are illustrated in Figure \ref{fig:distribution_plots_sparse}. These results demonstrate that our theory remains valid in the sparse edge setting.

\section{Analysis of The Statistician Coauthorship Network}
\label{sec:real_data}

We demonstrate the usefulness of our statistical inference results via an application to a statistician coauthorship network. The data is collected by the authors of \cite{ji2016coauthorship}, consisting of authorship information of papers published in four major statistical journals from 2003 to 2012. In this coauthorship network, each node represents a statistician, and two statisticians are connected if they have coauthored at least 2 papers. Our analysis focuses on the largest connected component, with $n = 236$ statisticians and $296$ edges. Consistently with \cite{ji2016coauthorship}, we use the scree-plot and choose to fit a two-dimensional latent space model on the coauthorship network. 

\begin{figure}[ht]
    \centering
    \includegraphics[width = 0.8 \linewidth]{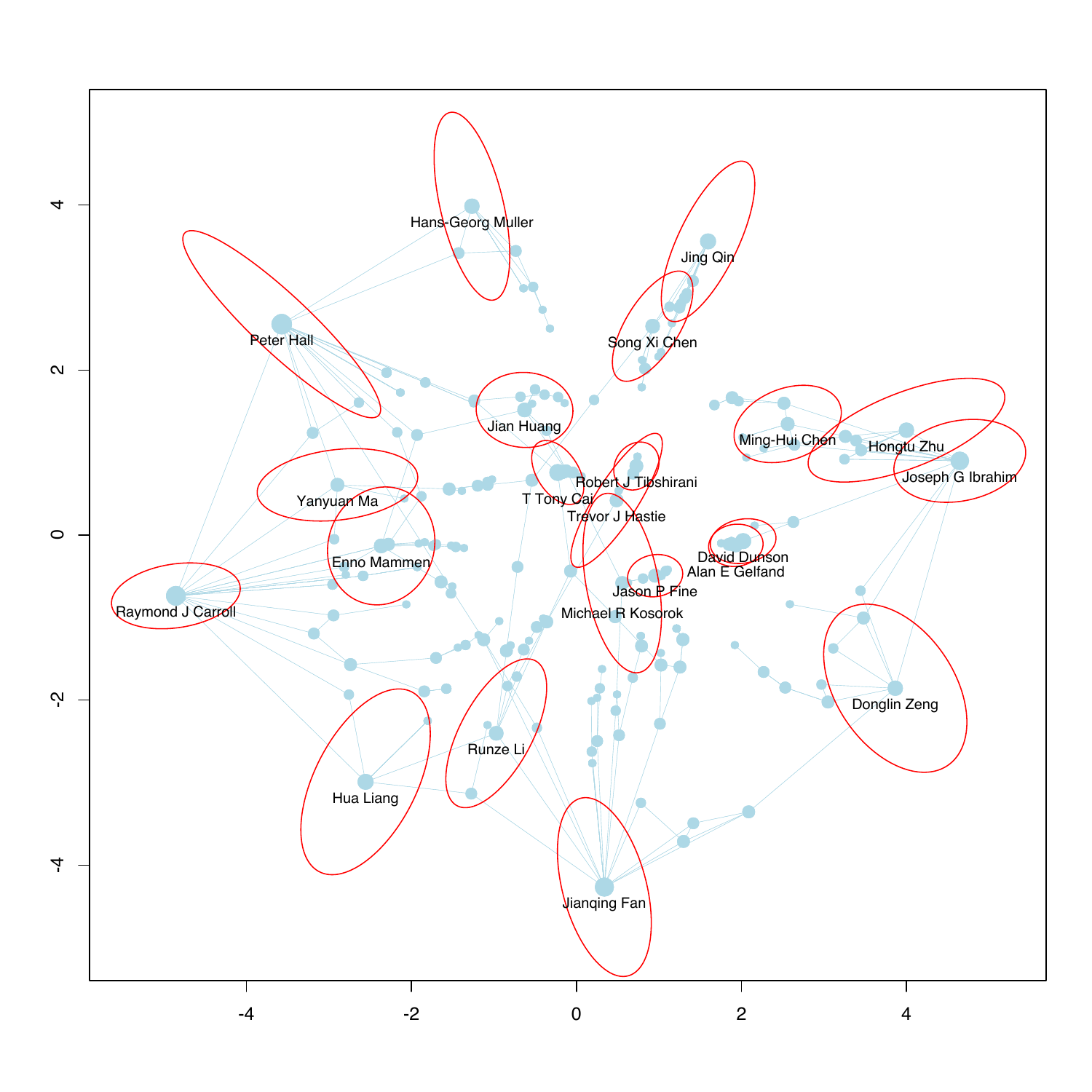}
    \caption{The estimated latent positions with 95\% confidence regions of hub nodes. Nodes are scaled by degrees.}
    \label{fig:estimated_Z}
\end{figure}

The estimated latent positions are visualized in Figure \ref{fig:estimated_Z}, providing insights into the collaboration patterns and research interests of statisticians within the network. The node sizes are proportional to the node degrees, with larger nodes representing higher degrees. We highlight those hub nodes (node degrees greater than 5) with the corresponding statisticians' names for better interpretation. A direct examination of the point estimates shows that the estimated two-dimensional latent positions uncover certain community structures among statisticians, which include groups such as non-parametric and semi-parametric (left), high-dimensional statistics (bottom), Bayesian statistics and Biostatistics (right), and  
statistical learning (middle). These observations align with the findings from related studies \citep{ji2016coauthorship}. However, our inference theories enable us to investigate beyond point estimation. We can construct confidence regions for the estimated latent positions, offering deeper insights into the network's nuances. In Figure \ref{fig:estimated_Z}, the 95\% confidence regions for the latent positions of the hub nodes are indicated by the red elliptical contours.

\begin{figure}[ht]
    \centering
    \includegraphics[width = 0.8 \linewidth]{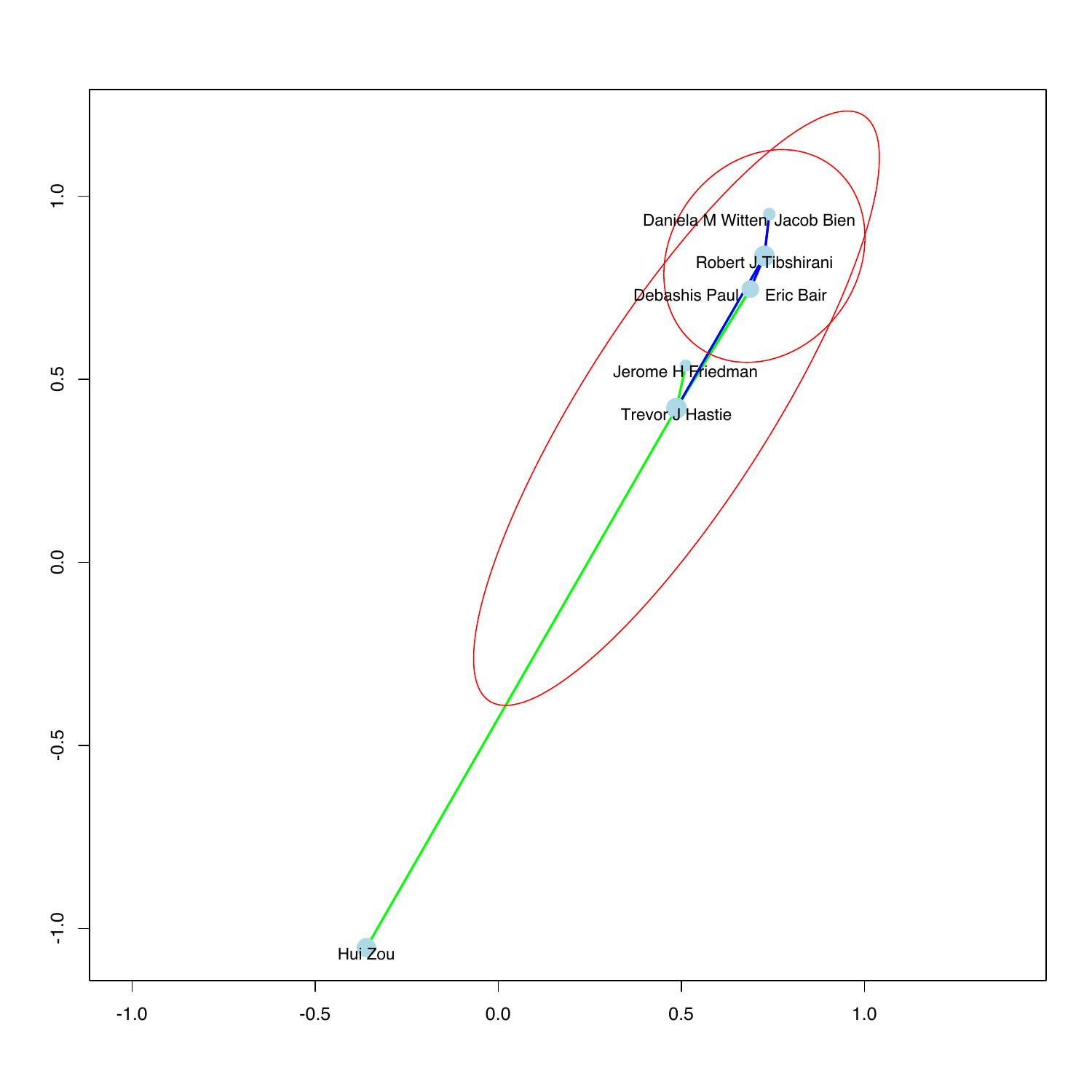}
    \caption{An example involving two statisticians, Robert J Tibshirani and Trevor J Hastie, with their collaborators. The blue and green lines are collaborations of Trevor J Hastie and Robert J Tibshirani, respectively.  The 95\% confidence regions of the two statisticians are plotted in red.}
    \label{fig:case_study}
\end{figure}

We further provide a statistical interpretation of the elliptical confidence regions regarding their sizes, shapes, and orientations. Based on Equations \eqref{eq:variance_estimators} and \eqref{eq:derivatives}, the key terms influencing the variance of the latent position estimator are those contributed by their neighboring nodes. In particular, the size of an ellipse decreases with an increase in the node degree but expands with a greater variance of the latent positions of the neighboring nodes. Thus, the estimation uncertainty of a statistician's latent position decreases with the number of collaborators while increasing with greater cross-domain collaborations. Moreover, from Equations \eqref{eq:variance_estimators} and \eqref{eq:derivatives}, the orientation of an ellipse's major and minor axes is largely influenced by the latent positions of collaborators. Additionally, the degree of overlap between ellipses could shed light on the similarity in collaboration patterns or research interests between statisticians.  

To further illustrate this, Figure \ref{fig:case_study} shows an example involving two statisticians, Dr. Robert J. Tibshirani and Dr. Trevor J. Hastie, along with their collaborators in the considered dataset.
In this network, although these two researchers' nodes share the same degrees and have latent positions relatively close to each other, their ellipses exhibit significant differences in sizes and orientations. This dissimilarity is primarily influenced by the estimated latent positions of their collaborators in the network, as indicated in the figure.

\section{Discussion}

In this paper, we address the crucial statistical inference problems of the latent space models for network data. Adopting a flexible analysis framework utilizing the Lagrange-adjusted Hessian matrix, which has not been introduced in the existing network literature, we prove the first uniform consistency and asymptotic distribution results for maximum likelihood estimators in a broad class of latent space models, accommodating different edge types and link functions. Furthermore, extensions have been established for two realistic yet challenging scenarios concerning edge dependency and sparsity. We conduct extensive simulation studies to validate the theoretical results and provide a data application focused on a statistician coauthorship network, demonstrating how the established theories can provide valuable insights into network structure beyond point estimation.


Our analysis techniques point to several promising future directions for downstream inference problems beyond link prediction. One direction involves network testing and node testing problems. Multi-network comparisons have been recently studied in the RDPG literature \citep{athreya2017statistical}, and two sample tests of a set of nodes have been studied in neuroimaging applications \citep{li2018spatially}. Our theoretical results can potentially be utilized to develop similar tests for a wide range of network data. {
Another direction concerns complex network data with node-level and/or edge-level covariates. For example, our theoretical results can be employed to conduct statistical inference for network-assisted prediction by accounting for the estimation uncertainty of latent positions, where the estimated network latent positions are used alongside node covariates to facilitate prediction \citep{lunde2023conformal}. 
Moreover, our proof techniques would provide useful analytic tools for conducting inference in the joint modeling of network data with node covariates \citep{zhang2022joint, li2025high} and edge covariates \citep{ma2020universal}, which are of great future interest.
}

%


\begin{supplement}
\stitle{Supplement to ``Statistical Inference on Latent Space Models for Network
Data''}  
\sdescription{The Supplementary Material contains additional simulation results and the proofs of main results.}
\end{supplement}




\bibliographystyle{imsart-nameyear}
\bibliography{ref}

%
%

\end{document}